\documentclass[12pt]{amsart}
\usepackage{amsthm}
\usepackage{amsmath, mathrsfs}
\usepackage{color}
\usepackage{amssymb}
\usepackage{hyperref}
\usepackage{enumerate}
\usepackage{latexsym,array}
\usepackage{amsfonts}
\usepackage{shadow}

\newtheorem{Pa}{Paper}[section]
\newtheorem{Tm}[Pa]{{\bf Theorem}}
\newtheorem{La}[Pa]{{\bf Lemma}}
\newtheorem{Dn}[Pa]{{\bf Definition}}
\newtheorem{Cy}[Pa]{{\bf Corollary}}

\newtheorem{Rk}[Pa]{{\bf Remark}}
\newtheorem{Pn}[Pa]{{\bf Proposition}}
\newtheorem{Pb}[Pa]{{\bf Problem}}

\newtheorem{Ex}[Pa]{{\bf Example}}

\newcommand{\w}{\omega}

\def\C{\mathbb C}
\author[D. Alpay]{Daniel Alpay}
\address{(DA) Department of Mathematics\\
Ben–Gurion University of the Negev\\
Beer-Sheva 84105 Israel} \email{dany@math.bgu.ac.il}

\author[F. Colombo]{Fabrizio Colombo}
\address{(FC) Politecnico di
Milano\\Dipartimento di Matematica\\Via E. Bonardi, 9\\20133
Milano, Italy}
\email{fabrizio.colombo@polimi.it}
\author[I. Sabadini]{Irene Sabadini}
\address{(IS) Politecnico di
Milano\\Dipartimento di Matematica\\Via E. Bonardi, 9\\20133
Milano, Italy}
\email{irene.sabadini@polimi.it}
\title[Pontryagin de Branges Rovnyak spaces]
{Pontryagin de Branges Rovnyak spaces of slice hyperholomorphic
functions} \oddsidemargin 0.2in \evensidemargin 0.2in \topmargin
-0.5in \textwidth 15.5truecm \textheight 23truecm

\def\R{\mathbb R}

\def\C{\mathbb C}

\def\(s){\mathscr S(\R\times\R)}



\thanks{The authors are grateful to Prof. V. Bolotnikov for the useful discussions and for a careful reading of the manuscript. D. Alpay thanks the Earl Katz family for endowing the chair
which supported his research, and the Binational Science
Foundation Grant number 2010117. F. Colombo and I. Sabadini acknowledge the Center for Advanced Studies of the Mathematical Department of the Ben-Gurion University of the Negev for the support and the kind hospitality during the period in which this paper has been written.}
\begin{document}
\maketitle \tableofcontents
\parindent 0cm
\begin{abstract}
We study reproducing kernel Hilbert  and Pontryagin spaces of slice
hyperholomorphic functions which are analogs of the Hilbert
spaces of analytic functions introduced by de Branges and Rovnyak.
In the first part of the paper we focus on the case of Hilbert
spaces, and introduce in particular a version of the Hardy space. Then we
define Blaschke factors and Blaschke products and we consider an
interpolation problem. In the second part of
the paper we turn to the case of Pontryagin spaces. We first
 prove some results from the
theory of Pontryagin spaces in the quaternionic setting and, in particular, a theorem of
Shmulyan on densely defined contractive linear relations. We then study realizations of generalized
Schur functions and of generalized Carath\'eodory functions.
\end{abstract}
\section{Introduction}
\setcounter{equation}{0}
Functions $s$ analytic in the open unit disk $\mathbb D$ and contractive there,
or equivalently such that the kernel
\[
\frac{1-s(z)s(w)^*}{1-zw^*}
\]
is positive definite in $\mathbb D$,
play an important role in operator theory, and their study is a
part of a field called Schur analysis.
The present work is a continuation of \cite{acs1}, and deals with various
aspects of Schur analysis in the case of slice hyperholomorphic functions.
To review the classical case, and to present the outline of the paper, we first
recall a definition:
A signature matrix is a matrix $J$ (say with complex entries; in the sequel
quaternionic entries will be allowed) which
is both self-adjoint and unitary. We denote by ${\rm sq}_-  J$ the
multiplicity (possibly equal to $0$) of the eigenvalue $-1$.
Let now $J_1$ and $J_2$ be two signature matrices, belonging to
$\mathbb C^{N\times N}$ and $\mathbb C^{M\times M}$ respectively, and
assume that
\[
{\rm sq}_- J_1={\rm sq}_-J_2.
\]
Functions $\Theta$ which are $\mathbb C^{M\times N}$-valued and
meromorphic in $\mathbb D$, and such that the kernel
\begin{equation}
\label{eq:kerneltheta}
K_\Theta(z,w)=\frac{J_2-\Theta(z)J_1\Theta(w)^*}{1-zw^*}
\end{equation}
has a finite number of negative squares in $\mathbb D$ are called generalized Schur functions, and have been studied by Krein
and Langer in a long series of papers; see for instance
\cite{kl2,MR518342,MR614775,MR563344,kl3}. These authors consider
also the case of operator-valued functions and other classes, in
particular, kernels of the form
\begin{equation}
\label{eq:kerneltheta1}
k_\varphi(z,w)=\frac{\varphi(z)J+J\varphi(w)^*}{1-zw^*},
\end{equation}
where $\varphi$ is $\mathbb C^{N\times N}$ valued and analytic in
a neighborhood of the origin, and $J\in\mathbb C^{N\times N}$ is a
signature matrix, and the counterparts of these kernels when the
open unit disk is replaced by the open upper half-plane.
Meromorphic functions $\Theta$ for which the kernel
\eqref{eq:kerneltheta} has a finite number of negative squares
are called generalized Schur functions, and meromorphic functions
$\Theta$ for which the kernel \eqref{eq:kerneltheta1} has a finite
number of negative squares are called generalized Carath\'eodory
functions. Associated problems (such as realization
and interpolation questions) have been studied extensively.\\

As  mentioned above, a study of Schur analysis in the setting of
slice hyperholomorphic functions has been initiated recently
in \cite{acs1}, and it is the
purpose of the present paper to continue this study. The paper
\cite{acs1} was set in the Hilbert spaces framework, and presented
in particular the notions and properties of Schur multipliers, de
Branges Rovnyak space, and coisometric realizations in the slice
hyperholomorphic setting. In the first part of this work we also focus
on the Hilbert space case, while in the second part we consider
the case of indefinite inner product spaces.\\

To set the present work into perspective we recall that the theory
of slice hyperholomorphic functions represents a novelty with
respect to other theories of hyperholomorphic functions that can
be defined in the quaternionic setting since it allows the
definition the quaternionic functional calculus and its
associated $S$-resolvent operator. The importance of the
$S$-resolvent operator, in the context of this paper, is the
definition of the quaternionic version of the operator
$(I-zA)^{-1}$ that appears in the realization function
$s(z)=D+zC(I-zA)^{-1}B$. It turns out that when $A$ is a
quaternionic matrix and $p$ is a quaternion then $(I-pA)^{-1}$
has to be replaced by
$$
(I-pA)^{-\star}=( I  -\bar p A)(|p|^2A^2-2{\rm Re}(p) A+I  )^{-1}
$$
 which is equal to $p^{-1}S^{-1}_R(p^{-1},A)$ where
$S^{-1}_R(p^{-1},A)$ is the right $S$-resolvent operator
associated to the quaternionic matrix $A$.
\\

Slice hyperholomorphic functions have two main formulations
according to the fact that the functions we consider are defined on quaternions and are
quaternion-valued, in this case the functions are called slice regular, see \cite{MR2353257, MR2555912, MR2742644} or the  functions are defined on the
Euclidean space $\mathbb{R}^{N+1}$ and have values in the Clifford
Algebra $\mathbb{R}_N$ and are also called slice monogenic functions, see \cite{MR2520116, MR2684426}. We also point out that there exists a non constant
coefficients differential operator whose kernel contains slice
hyperholomorphic
functions defined on suitable domains, see \cite{GLOBAL}.\\

 Slice hyperholomorphicity has applications in operator theory: specifically, in the case of quaternions, it allows the definition  of a quaternionic functional calculus, see e.g. \cite{MR2735309, MR2496568,  MR2803786}; while
slice monogenic functions admit a functional calculus for
$n$-tuples of operators,  see
\cite{MR2402108, MR2720712, SCFUNCTIONAL}. The book
\cite{MR2752913} collects some of the main results on the theory
of slice hyperholomorphic functions and the related functional
calculi.\\

 Finally we mention the paper
\cite{MR2124899,MR2240272,asv-cras}, where Schur multipliers were
introduced and studied in the quaternionic setting using the
Cauchy-Kovalesvkaya product and series of Fueter polynomials, and
the papers \cite{MR733953,MR2205693,MR2131920}, which treat
various aspects of a theory of linear systems in the quaternionic
setting. Our approach is quite different from the methods used
there.\\

The paper consists of nine sections besides the introduction, and
its outline is as follows: In Sections 2  and 3 we review some
basic definitions on slice hyperholomorphic functions. In Section
4 we discuss the notion of multipliers in the case of reproducing
kernel Hilbert spaces of slice hyperholomorphic functions. In
Section 5 we discuss the Hardy space in the present setting, and
introduce Blaschke products. Interpolation in the Hardy space is
studied in Section 6. Sections 7-11 are in the setting of
indefinite metric spaces. A number of facts on quaternionic
Pointryagin spaces as well as a proof of a theorem of Shmulyan on
relations are proved in Section 7. Negative squares are discussed
in Section 8, while Section 9 introduces generalized Schur
functions and discusses their realizations. We also consider in
this section the finite dimensional case. Finally, we briefly
discuss in Section 10 the case of generalized Carath\'eodory
functions.

\section{Slice hyperholomorphic functions}
In the literature there are several notions of quaternion valued
hyperholomorphic functions. In this paper we consider a notion
which includes power series in the quaternionic variable, the
so-called slice regular or slice hyperholomorphic functions, see
\cite{MR2752913}. In order to introduce the class of slice
hyperholomorphic functions, we  fix some preliminary notations.
By $\mathbb{H}$ we denote the algebra of real quaternions
$p=x_0+ix_1+jx_2+kx_3$. A quaternion can also be written as
$p={\rm Re}(p)+{\rm Im}(p)$ where $x_0={\rm Re}(p)$ and
$ix_1+jx_2+kx_3={\rm Im}(p)$ but also as $q={\rm Re}(p)+I_p|{\rm
Im}(p)|$ where $I_p={\rm Im}(p)/|{\rm Im}(p)|$, as long as $p$ is
non real. The element $I$ belongs to the 2-sphere
$$
\mathbb{S}=\{p=x_1i+x_2j+x_3k\ :\ x_1^2+x_2^2+x_3^2=1\}
$$
of unit purely imaginary quaternions.
\begin{Dn}
Let $\Omega\subseteq\mathbb{H}$ be an open set and let $f:\
\Omega\to\mathbb{H}$ be a real differentiable function. Let
$I\in\mathbb{S}$ and let $f_I$ be the restriction of $f$ to the
complex plane $\mathbb{C}_I := \mathbb{R}+I\mathbb{R}$ passing
through $1$ and $I$ and denote by $x+Iy$ an element in
$\mathbb{C}_I$.
\begin{enumerate}
\item
 We say that $f$ is a left slice regular function
(or slice regular or slice hyperholomorphic)  if, for every
$I\in\mathbb{S}$, we have:
$$
\frac{1}{2}\left(\frac{\partial }{\partial x}+I\frac{\partial
}{\partial y}\right)f_I(x+Iy)=0.
$$
\item
We say that $f$ is right slice regular function (or right slice
hyperholomorphic) if, for every $I\in\mathbb{S}$, we have
$$
\frac{1}{2}\left(\frac{\partial }{\partial
x}f_I(x+Iy)+\frac{\partial }{\partial y}f_I(x+Iy)I\right)=0.
$$
\end{enumerate}
\end{Dn}
\begin{Dn}
The set of all elements of the form ${\rm Re}(p)+J |{\rm Im}(p)|$
when $J$ varies in $\mathbb{S}$ is denoted by $[p]$ and is called
 the 2-sphere defined by $p$.
\end{Dn}
The most important feature of slice hyperholomorphic functions is
that, on a suitable class of open sets described below, they can
be reconstructed by knowing their values on a complex plane
$\mathbb{C}_I$ by the so-called Representation Formula.
\begin{Dn}
Let $\Omega$ be a domain in $\mathbb{H}$. We say that $\Omega$ is
a \textnormal{slice domain} (s-domain for short) if $\Omega \cap
\mathbb{R}$ is non empty and if $\Omega\cap \mathbb{C}_I$ is a
domain in $\mathbb{C}_I$ for all $I \in \mathbb{S}$. We say that
$\Omega$ is \textnormal{axially symmetric} if, for all $p\in
\Omega$, the 2-sphere $[p]$ is contained in $\Omega$.
\end{Dn}

\begin{Tm}[Representation Formula]\label{formula}
Let $\Omega\subseteq \mathbb{H}$ be an axially symmetric s-domain.
Let $f$ be a left slice regular function on  $\Omega\subseteq
\mathbb{H}$.
 Then the following equality holds for all $p=x+I_p y \in \Omega$:
\begin{equation}\label{strutturaquat}
f(p)=f(x+I_p y) =\frac{1}{2}\Big[   f(z)+f(\overline{z})\Big]
+\frac{1}{2}I_pI\Big[ f(\overline{z})-f(z)\Big],
\end{equation}
where $z:=x+Iy$, $\overline{z}:=x-Iy\in\Omega\cap\mathbb{C}_I$.
Let $f$ be a right slice regular function on  $\Omega\subseteq
\mathbb{H}$. Then the following equality holds for all $p=x+I_p y
\in \Omega$:
\begin{equation}\label{distribution}
f(x+I_p y) =\frac{1}{2}\Big[   f(z)+f(\overline{z})\Big]
+\frac{1}{2}\Big[ f(\overline{z})-f(z)\Big]II_p.
\end{equation}
\end{Tm}
The Representation Formula allows to extend any function $f: \
\Omega\subseteq\mathbb{C}_I\to\mathbb{H}$ defined on an axially
symmetric s-domain  $\Omega$  and in
the kernel of the corresponding Cauchy-Riemann operator to a
function $f: \
\widetilde{\Omega}\subseteq\mathbb{H}\to\mathbb{H}$ slice
hyperholomorphic where $\widetilde{\Omega}$ is the smallest
axially symmetric open set in $\mathbb{H}$ containing $\Omega$.  Using the above notations, the extension is obtained by
means of the {\em extension operator}
\begin{equation}\label{ext}
{\rm ext}(f)(p):= \frac{1}{2}\Big[   f(z)+f(\overline{z})\Big]
+\frac{1}{2}I_pI\Big[ f(\overline{z})-f(z)\Big],\quad z,\bar
z\in\Omega\cap\mathbb{C}_I, \ p\in\widetilde\Omega.
\end{equation}
For example, in the case of the kernel associated to the Hardy
space, the extension operator applied to the function
 $\sum_{n=0}^\infty z^n\bar w^n$ gives (see Proposition 5.3 in \cite{acs1}):
\begin{Pn}
Let $p$ and $q$ be quaternionic variables.
The sum of the series $\sum_{n=0}^{+\infty} p^n\bar q^n$ is the
function $k(p,q)$ given by
\begin{equation}\label{kernel}
k(p,q)=(1-2{\rm Re}(q) p+|q|^2p^2)^{-1}(1-pq)=(1-\bar p \bar q)
(1-2{\rm Re}(p) \bar q+|p|^2\bar q^2)^{-1}.
\end{equation}
The kernel $k(p,q)$ is defined for all $p$ outside the 2-sphere
defined by $[q^{-1}]$ (or, equivalently, for all $q$ outside the
2-sphere $[p^{-1}]$. Moreover:
\begin{itemize}
\item[a)] $k(p,q)$ is slice hyperholomorphic in $p$ and right slice hyperholomorphic in $\bar q$;
\item[b)] $\overline{k(p,q)}=k(q,p)$.
\end{itemize}
\end{Pn}

The function $k(p,q)$ in the preceding proposition is positive definite, and is
the reproducing kernel of the slice hyperholomorphic counterpart of the Hardy space
$\mathbf H_2(\mathbb{B})$ of functions analytic in the open unit ball $\mathbb{B}$, see \cite{acs1} and
Section \ref{secHardy} below.
\begin{Rk}\label{starlandr} {\rm  The two possible expressions for $k(p,q)$ given in (\ref{kernel}) correspond to the left slice regular reciprocal of $1-p\bar q$ in the variable $p$ and to the right slice regular reciprocal in the variable $\bar q$, (see the discussion in \cite[Proposition 5.3]{acs1}) and these two reciprocals coincide. Thus, in the sequel, we will often write $(1-p\bar q)^{-\star}$ instead of $k(p,q)$.
}
\end{Rk}
\begin{Rk}{\rm
Note that whenever a function $k(p,q)$ is slice hyperholomorphic
in $p$ and is Hermitian, then $k(p,q)$
is right slice hyperholomorphic in $\bar q$.}
\end{Rk}

\section{Slice hyperholomorphic multiplication}
\setcounter{equation}{0} We recall that, given two left slice
hyperholomorphic functions $f$, $g$, it is possible to introduce
a binary operation called the $\star$-product, such that $f\star g$ is a slice hyperholomorphic function.
Similarly, given two right slice hyperholomorphic functions, we
can define their $\star$-product. When considering in same
formula both the products, it may be useful to distinguish
between them and in this case we will write $\star_l$ or
$\star_r$ according to the fact that we are using the left or the
right slice regular product. When there is no subscript, we will
mean that we are considering the left $\star$-product.\\

Let $f,g:\ \Omega \subseteq\mathbb{H}$ be slice regular functions such that
their restrictions to the complex plane $\mathbb{C}_I$ can be written as
$f_I(z)=F(z)+G(z)J$,
$g_I(z)=H(z)+L(z)J$ where $J\in\mathbb{S}$, $J\perp I$. The functions $F$,
$G$, $H$, $L$ are holomorphic functions of the variable $z\in
\Omega \cap \mathbb{C}_I$ and they exist by the splitting lemma, see \cite{MR2752913}, p. 117.
The $\star_l$-product of  $f$ and $g$ is defined as the unique
left slice hyperholomorphic function whose restriction to the
complex plane $\mathbb{C}_I$ is given by
\begin{equation}\label{starproduct}
(F(z)+G(z)J)\star_l(H(z)+L(z)J):=
(F(z)H(z)-G(z)\overline{L(\bar z)})+(G(z)\overline{H(\bar z)}+F(z)L(z))J.
\end{equation}

Pointwise
multiplication and slice multiplication are different, but they can be related as in the following result, \cite[Proposition
4.3.22]{MR2752913}:
\begin{Pn}\label{pointwiseproduct}
Let $U \subseteq \mathbb{H}$ be an axially symmetric s-domain,
 $f, g : U \to \mathbb{H}$ be slice hyperholomorphic functions
and let us assume that $f(p)\not=0$. Then
\begin{equation}
\label{productfg}
(f \star g)(p) = f(p) g(f(p)^{-1}pf(p)),
\end{equation}
for all $p\in U$.
\end{Pn}
\begin{Rk}{\rm
The transformation $p\to f(p)^{-1}pf(p)$ is clearly a rotation in $\mathbb{H}$,
since $|p|=|f(p)^{-1}pf(p)|$
and allows to rewrite the $\star$-product as a pointwise product.\\
Note also that if $f\star g(p)=0$ then either $f(p)=0$ or $g(f(p)^{-1}pf(p))=0$.
}
\end{Rk}

As a consequence of Proposition \ref{pointwiseproduct} one has:

\begin{Cy}\label{rotation}
If
$
\lim_{r\to 1}|f(re^{I\theta})|=1,
$
for all $I$ fixed in $\mathbb{S}$, then
$$
\lim_{r\to 1} |f\star g(re^{I\theta})|=| g(e^{I'\theta})|,
$$
where $\theta\in[0,2\pi)$,  and $I'\in\mathbb{S}$ depends on
$\theta$ and $f$.
\end{Cy}
\begin{proof}
Set $b=f(re^{I\theta})$.
We can write $b=Re^{J\alpha}$ for suitable $R, J,\alpha$ and, by hypothesis, we can assume that $b\not=0$ when $r\to 1$, thus $b^{-1}$ exists.
We have
\[
\begin{split}
&b^{-1}re^{I\theta}b=e^{-J\alpha}(r e^{I\theta})e^{J\alpha}=r(\cos\alpha -J\sin\alpha)(\cos\theta +I\sin\theta)(\cos\alpha + J\sin\alpha)\\
&=r(\cos\theta +I\cos^2\alpha \sin\theta -JI \cos\alpha \sin\alpha \sin\theta  + IJ \cos\alpha \sin\alpha \sin\theta-JIJ  \sin^2\alpha \sin\theta)\\
&=r(\cos\theta +\cos\alpha e^{-J\alpha}I\sin\theta+e^{-J\alpha}IJ\sin\alpha\sin\theta)\\
&=r( \cos\theta + e^{-J\alpha}I e^{J\alpha}\sin\theta)= r(\cos\theta + I'\sin\theta),\\
\end{split}
\]
where $I'=e^{-J\alpha}I e^{J\alpha}$.
 Then, the result immediately follows
from the equalities:
$$
\lim_{r\to 1} |f\star g(re^{I\theta})|= \lim_{r\to 1}
|f(re^{I\theta}) g(b^{-1}re^{I\theta}b)| =\lim_{r\to 1} |
g(re^{I'\theta})|=| g(e^{I'\theta})|.
$$
\end{proof}

Given a left slice regular function $f$ it is possible to construct its slice regular
reciprocal, which is denoted by $f^{-\star}$. The
general construction can be found in \cite{MR2752913}. In this paper we will be in need
of the reciprocal of a polynomial
or a power series with center at the origin that can be described in the easier way illustrated below.
\begin{Dn}
Given $f(p)=\sum_{n=0}^\infty p^n a_n$, let
us set
$$
f^c(p)=\sum_{n=0}^\infty p^n \bar a_n,\qquad  f^s(p)=(f^c\star f)(p
)=\sum_{n=0}^\infty p^nc_n,\quad
c_n=\sum_{r=0}^n a_r\bar a_{n-r},
$$
where the series converge.
The left slice hyperholomorphic reciprocal of $f$
is then defined as
$$
f^{-\star}:=(f^s)^{-1}f^c.
$$
\end{Dn}
 In an analogous way one can define
the right slice hyperholomorphic reciprocal $
f^{-\star}:=f^c(f^s)^{-1},
$ of a right slice regular function  $f(q)=\sum_n a_n q^n $.
Note that the series $f^s$ has real coefficients.
\begin{Rk}{\rm
Let $\Omega$ be an axially symmetric open set. We recall that if $f$ is
left slice hyperholomorphic in $q\in \Omega$ then $\overline{f(q)}$ is
right slice hyperholomorphic in $\overline{q}$. This fact follows
immediately from $(\partial_x+I\partial_y)f_I(x+Iy)=0$, since by
conjugation we get
$\overline{f_I(x+Iy)}(\partial_x-I\partial_y)=0$ for all $I\in
\mathbb{S}$. }
\end{Rk}
\begin{La}\label{starmult}
Let $\Omega$ be an axially symmetric s-domain and let $f,g:\Omega
\to \mathbb{H}$ be two left slice hyperholomorphic functions.
Then $$\overline{f\star_l g}=\overline{g}\star_r \overline{f},$$
where $\star_l$, $\star_r$ are the left and right $\star$-products with respect to $q$ and $\bar q$, respectively.
\end{La}
\begin{proof}
Let $f_I(z)=F(z)+G(z)J$,
$g_I(z)=H(z)+L(z)J$ be the restrictions of $f$ and $g$ to the
complex plane $\mathbb{C}_I$, respectively. The functions $F$,
$G$, $H$, $L$ are holomorphic functions of the variable $z\in
\Omega \cap \mathbb{C}_I$ which exist by the splitting lemma and $J$ is an element in the sphere $\mathbb{S}$ orthogonal to $I$. The
$\star_r$-product of two right slice hyperholomorphic functions
$\overline{g}$ and $\overline{f}$ in the variable $\overline{q}$
is defined as the unique right slice hyperholomorphic function
whose restriction to a complex plane $\mathbb{C}_I$ is given by
$$
(\overline{H(z)}-J \ \overline{L(z)})\star_r(\overline{F(z)}-J\
\overline{G(z)}):= (\overline{H(z)}\
\overline{F(z)}-L(\bar z)\overline{G(z)})-J(\overline{L(z)}\
\overline{F(z)}+H(\bar z)\overline{G(z)}).
$$
Thus, comparing with (\ref{starproduct}), it is clear that
$$
\overline{f_I\star_l g_I}=\overline{g_I}\star_r \overline{f_I},
$$
and the statement follows by taking the unique right slice
hyperholomorphic extension.
\end{proof}
\begin{Rk}{\rm
For the sake of completeness, we adapt some of the previous definitions in the case we consider matrix valued functions.
We will say that a real differentiable function $f:\Omega\subseteq\mathbb{H}\to \mathbb{H}^{N\times M}$ is left (resp. right)  slice hyperholomorphic if
and only if for any linear and continuous functional $\Lambda$ acting on   $\mathbb{H}^{N\times M}$, the function $\Lambda f$ is left (resp. right) slice hyperholomorphic in $\Omega$.
If, in particular, $\Omega=\mathbb{B}$, then it can be shown with standard techniques that  $f$ is left slice hyperholomorphic if and only if
$f(p)=\sum_{n=0}^\infty p^n A_n$, where $A_n\in\mathbb{H}^{N\times M}$ and the series converges in $\mathbb{B}$. Let $f:\ \mathbb{B}\to\mathbb{H}^{N\times M}$,
 $g:\ \mathbb{B}\to\mathbb{H}^{M\times L}$ be left slice hyperholomorphic and let $f(p)=\sum_{n=0}^\infty p^n A_n$,
 $g(p)=\sum_{n=0}^\infty p^n B_n$. The $\star$-product of  $f$ and  $g$  is defined as $f\star g:=\sum_{n=0}^\infty p^n C_n$ where $C_n=\sum_{r=0}^n A_rB_{n-r}$.  Analogous definitions can be given in the case we consider right slice hyperholomorphic functions.}
\end{Rk}
\begin{Rk}{\rm
When considering the function $\sum_{n=0}^\infty p^nA^n$ where $A\in\mathbb{H}^{N\times N}$ and $|p|<1/\|A\|$,  or, more in general, $A$ is a bounded right linear quaternionic operator from a quaternionic Hilbert space to itself, then $(I-pA)^{-\star}$ can be constructed using the functional calculus (see \cite[Proposition 2.16]{acs1}): it is sufficient to construct the right slice regular inverse of $1-pq$ with respect to $q$ and then substitute $q$ by the operator $A$. Note that we write $(I-pA)^{-\star}$ using the symbol $\star$ instead of $\star_r$ for simplicity and the discussion in Remark \ref{starlandr} justifies this abuse of notation.
}
\end{Rk}

\section{Multipliers in reproducing kernel Hilbert spaces}\label{positivity}
In this section we study the multiplication operators and their adjoints, we show that positivity implies the slice hyperholomorphicity for
a class of functions and, finally, we prove that if a kernel is positive and slice hyperholomorphic then the corresponding reproducing kernel Hilbert space consists of slice hyperholomorphic functions.

\setcounter{equation}{0} Let us begin by recalling the following definition,
see \cite{as3}:
\begin{Dn}\label{rkHs}
A quaternionic Hilbert space $\mathscr{H}$ of
$\mathbb{H}^{N}$-valued functions defined on an open set
$\Omega\subseteq\mathbb{H}$ is called a reproducing kernel
quaternionic Hilbert space if there exists a $\mathbb{H}^{N\times
N}$-valued function defined on $\Omega\times\Omega$ such that:
\begin{enumerate}
\item  For every $q\in\Omega$ and $a\in\mathbb{H}^{N}$ the function
$
p\mapsto K(p,q)a
$
belongs to $\mathscr{H}$.
\item  For every $f\in\mathscr{H}$, $q\in\Omega$ and $a\in\mathbb{H}^{N}$
$$
\langle f, K(\cdot , q) a\rangle_{\mathscr H} =a^* f(q).
$$
\end{enumerate}
\end{Dn}
The function $K(p,q)$ is called the reproducing kernel of the
space. As observed in \cite{as3}, Definition \ref{rkHs}, one may
ask the weaker requirement  that $\mathscr{H}$ is a quaternionic
pre-Hilbert space. However, the next result proven in \cite{as3},
guarantees that a reproducing kernel quaternionic pre-Hilbert
space has a unique completion as a reproducing kernel
quaternionic Hilbert space, which will be denoted by
$\mathscr{H}(K)$.
\begin{Tm}
Given an $\mathbb{H}^{N\times N}$-valued function $K(p,q)$
positive on a set $\Omega \subset  \mathbb{H}$, there exists a
uniquely defined reproducing kernel quaternionic Hilbert space of
$\mathbb{H}^{N}$-valued function defined on  $\Omega$ and with
reproducing kernel $K(p,q)$.
\end{Tm}
Let us recall that $\mathscr{H}(K)$ is the completion of the
linear span $\stackrel{\circ}{\mathscr{H}(K)}$ of functions of
the form
\begin{equation}\label{span}
p\mapsto K(p,q)a, \qquad q\in\Omega,\ a\in\mathbb{H}^N,
\end{equation}
with the inner product
\begin{equation}\label{inner}
\langle K(\cdot, q)a, K(\cdot,
s)b\rangle_{\stackrel{\circ}{\mathscr{H}(K)}}:= b^*K(s,q)a.
\end{equation}

\begin{Pn}
Let $\phi$ be a slice hyperholomorphic function defined on an
axially symmetric s-domain $\Omega$ and with values in
$\mathbb{H}^{N\times M}$, and let $K_1(p,q)$ and $K_2(p,q)$ be
positive definite kernels in $\Omega$, respectively $\mathbb
H^{M\times M}$- and $\mathbb H^{N\times N}$-valued, and slice
hyperholomorphic in the variable $p$. Moreover,\\
$(1)$ Assume that the
slice multiplication operator
\[
M_\phi\,\,:\,\, \mathcal{H}(K_1)\to
\mathcal{H}(K_2)
\]
given by
\[
M_\phi\,\,:\,\, f\mapsto \phi\star f
\]
is continuous. Then, the adjoint operator is given by the formula:
$$
M^*_\phi(K_2(\cdot , q) d) =K_1(\cdot ,q )\star_r \phi^*(q) d.
$$
$(2)$ The multiplication operator $M_\phi$ is bounded and with
norm less or equal to $k$ if and only if the function
\begin{equation}\label{kernelK1K2star}
K_2(p,q)-\frac{1}{k^2}\phi(p)\star_l K_1(p,q)\star_r{\phi(q)^*}
\end{equation}
is positive on $\Omega$.
\end{Pn}
\begin{proof}
We compute the adjoint of the multiplication operator $M_\phi:\
\mathcal{H}(K_1)\to \mathcal{H}(K_2)$:
\[
\begin{split}
c^*(M^*_\phi(K_2(\cdot , q) d)(p) &= \langle
M_{\phi}^*(K_2(\cdot , q)d), K_1(\cdot ,p) c\rangle_{\mathcal{H}(K_1)}\\
&=\langle K_2(\cdot , q)d, \phi \star_l K_1(\cdot ,p) c\rangle_{\mathcal{H}(K_2)}\\
&= \langle \phi \star_l K_1(\cdot ,p) c, K_2(\cdot , q)d \rangle^*_{\mathcal{H}(K_2)}\\
&= (d^* (\phi(q)\star_l  K_1(q ,p))c)^*\\
&= c^*  (\phi(q)\star_l  K_1(q ,p))^* d.\\
\end{split}
\]
Now observe that by Lemma \ref{starmult} we have $(\phi(q)\star_l
K_1(q ,p))^*=  K_1(p,q )\star_r {\phi}^*(q)$ and so
$$
M^*_\phi(K_2(\cdot , q) d) =K_1( \cdot,q )\star_r \phi^*(q)d.
$$
The positivity of (\ref{kernelK1K2star}) follows from the positivity
of the operator $k^2-M_\phi M_\phi^*$. Conversely, if (\ref{kernelK1K2star}) is
positive, the standard argument shows that $\|M_\phi\|\le k$.
\end{proof}
\begin{Ex}
Let us consider the case in which the kernel $K$ is                                                                                               of the form
$$
K(p,q)=\sum_{n=0}^{\infty} p^n\overline{q}^n\alpha_n, \ \ \
\alpha_n\in \mathbb{R}, \ \ \forall n\in \mathbb{N}.
$$
Then we have
$$
\phi(p)\star_l K(p,q)=\sum_{n=0}^{\infty} p^n \phi(p)
\overline{q}^n\alpha_n
$$
and
$$
(\phi(p)\star_l  K(p,q))^*=\sum_{n=0}^{\infty} q^n \phi(p)^*
\overline{p}^n\alpha_n,
$$
from which we obtain
$$
\phi(q)\star_l (\phi(p)\star_l  K(p,q)) ^*= \phi(q)\star_l
\sum_{n=0}^{\infty} q^n \phi(p)^* \overline{p}^n\alpha_n
=\phi(q)\star_l K(q,p)\star_r \phi(p)^*.
$$
\end{Ex}

Recall that we defined in \cite{acs1} a Schur function to be a
function $S$ with values in $\mathbb
H^{N\times M}$, slice hyperholomorphic in $\mathbb B$ and such that
the kernel
\begin{equation}\label{ksh}
k_{S}(p,q)=\sum_{n=0}^\infty p^n(I_N-{
S}(p)S(q)^*)\overline{q}^n\\
= (I_N-S(p)S(q)^*) \star (1-p\overline{q})^{-\star}
\end{equation}
is positive on $\mathbb{B}$. We will show in Theorem \ref{2.6.3 Panorama} below that the converse, i.e. that positivity forces hyperholomorphicity, is
true for a subclass of slice hyperholomorphic functions. This
subclass is denoted by $\mathcal{N}$ and corresponds to those
functions $f$ such that $f: \mathbb{B}\cap \mathbb{C}_I\to
\mathbb{C}_I$ for any $I\in \mathbb{S}$. For these functions the pointwise multiplication of $f$ with a monomial of the form $p^n$ is
well defined and commutative since $f$ takes the complex plane $\mathbb{C}_{I_p}$ to itself and thus it behaves, on each plane, like a complex valued function.\\

We will be in need of the following preliminary result, see
\cite[Proposition 9.3]{as3}.
\begin{Pn}\label{2.3.9}
Let $K_1$ and $K_2$ be two positive functions on a set $\Omega$
with values in $\mathbb{H}^{N\times N}$ and $\mathbb{H}^{M\times
M}$, respectively. Let $\phi$ be a function defined on $\Omega$
and with values in $\mathbb{H}^{N\times M}$. The pointwise
multiplication operator by $\phi$ is bounded and with norm less or
equal to $k$ if and only if the function
\begin{equation}\label{kernelK1K2}
K_2(p,q)-\frac{1}{k^2}\phi(p) K_1(p,q){\phi(q)^*}
\end{equation}
is positive on $\Omega$.
\end{Pn}

\begin{Tm}\label{2.6.3 Panorama}
Let $S: \mathbb{B}\to \mathbb{H}^{N\times M}$ be a function such
that $S: \mathbb{B}\cap \mathbb{C}_I\to \mathbb{C}_I^{N\times M}$
for every $I\in \mathbb{S}$. The following are equivalent:
\begin{enumerate}
\item[(1)]
The function $\sum_{n=0}^\infty  p^n(I_N-S(p)S(q)^*)\bar q^n$ is
positive on $\mathbb{B}$.
\item[(2)]
The operator $M_S$ is a contraction from
$\mathbf{H}_2^{M}(\mathbb{B})$ to $\mathbf{H}_2^{N}(\mathbb{B})$.
\item[(3)]
$S$ is a Schur function belonging to $\mathcal{N}(\mathbb{B})$.
\end{enumerate}
\end{Tm}
\begin{proof}
The equivalence between (1) and (2) follows as in
\cite{MR2002b:47144} Theorem 2.6.3, and its proof is based on Proposition \ref{2.3.9}.
Indeed, let us set in (\ref{kernelK1K2})
$$
K_1(p,q)=I_M(1-p\bar q)^{-\star}, \qquad K_2(p,q)=I_N(1-p\bar
q)^{-\star}.
$$
We have:
$$
I_N(1-p\bar q)^{-\star}-\frac{1}{k^2}S(p)(1-p\bar q)^{-\star}
S(q)^*
$$
\begin{equation}\label{last}
=I_N\sum_{n=0}^\infty p^n \bar q^n
-\frac{1}{k^2}S(p)(\sum_{n=0}^\infty p^n \bar q^n ) S(q)^*;
\end{equation}
now observe that, by hypothesis, $S(p)$ commutes with $p^n$ since
$S$ takes the complex plane $\mathbb{C}_{I_p}$ to itself;
similarly, $S(q)^*$ commutes with $\bar q^n$. So we obtain that
(\ref{last}) is equal to:
$$
\sum_{n=0}^\infty p^n (I_N-\frac{1}{k^2}S(p) S(q)^*)\bar q^n.
$$
Thus, if (1) holds then by Proposition  \ref{2.3.9} we conclude
that $M_S$ is a contraction from $\mathbf{H}_2^{M}(\mathbb{B})$
to $\mathbf{H}_2^{N}(\mathbb{B})$. Conversely, if (2) holds, then
again Proposition \ref{2.3.9} allows to conclude that (1) holds.

The implication (3)$\Rightarrow$(2) follows from the fact that $S$ is a
Schur function. We show that (2) implies (3). The function $S$ is
slice hyperholomorphic since $Sc\in \mathbf{H}_2^N(\mathbb{B})$
for any $c\in \mathbb{H}^M$. Observe that the function $S$ is
contractive since $M^*_S$ acts as
$$
M^*_S\Big( (1-p\overline{q})^{-\star}d \Big)=
(1-p\overline{q})^{-\star}  S(q)^* d
$$
and it is a contraction.
\end{proof}

\begin{Dn}
A subset $\Omega$ of $\mathbb{B}$ is called a set of uniqueness if every slice
hyperholomorphic function on $\mathbb{B}$ which vanishes on
$\Omega$ is identically zero on $\mathbb{B}$.
\end{Dn}
\begin{Ex}
Any open subset $\Omega$ of $\mathbb{B}\cap \mathbb{C}_I$ is
a set of uniqueness. More in general, any subset $\Omega$ of
$\mathbb{B}\cap \mathbb{C}_I$ for $I\in \mathbb{S}$ having an
accumulation point in  $\mathbb{C}_I$ is a set of uniqueness.
\end{Ex}

\begin{Tm}
Let $\Omega$ be a set of uniqueness in $\mathbb{B}$ and let $S$
be a function defined on $\Omega$ such that $S: \Omega\cap
\mathbb{C}_I\to \mathbb{C}_I^{N\times M}$ for every $I\in
\mathbb{S}$. Then
$S$ can be extended slice hyperholomorphically to a Schur
function in $\mathcal{N}(\mathbb{B})$ if and only if
the kernel
 \begin{equation}\label{kernelstar}
 \sum_{n=0}^\infty  p^n(I_N-S(p)S(q)^*)\bar q^n
 \end{equation}
is positive on $\Omega$.
\end{Tm}
\begin{proof}
If $S$ can be extended hyperholomorphically to a Schur function, then the kernel (\ref{kernelstar}) is positive definite on $\Omega$.
We prove the converse.
Define the right linear quaternionic operator $T$  as
$$
T\Big( (1-p\overline{q})^{-\star}d
\Big)=(1-p\overline{q})^{-\star} S(q)^*d
$$
for $q\in \Omega$ and reason as in the proof of Theorem
\ref{2.6.3 Panorama}. By assumption the kernel
$\sum_{n=0}^\infty  p^n(I_N-S(p)S(q)^*)\bar q^n$ is positive thus
$T$ is well defined and contractive. Its domain is dense since
$\Omega$ is a set of uniqueness. So $T$ extends to a contraction
from $\mathbf{H}_2^M$ to $\mathbf{H}_2^N$. Its adjoint is a
contraction and for any $q\in \Omega$ and $F\in \mathbf{H}_2^N$
we have
\[
\begin{split}
\langle T^*F, (1-p\overline{q})^{-\star}d \rangle&=\langle F,
T\Big( (1-p\overline{q})^{-\star}d \Big)\rangle\\
&
=\langle F, (1-p\overline{q})^{-\star} S(q)^*d \rangle\\
& =d^* S(q)
F(q).
\end{split}
\]
Since we obtained a function equal to $S(q) F(q)$ on $\Omega$, the choice
$F=1$ shows that $S=T^*1$ is the restriction to $\Omega$ of a Schur
function.
\end{proof}

To conclude this section we show that if a kernel $K(p,q)$ is
positive and slice hyperholomorphic in $p$, then its
corresponding reproducing kernel Hilbert space consists of slice
hyperholomorphic functions.
\begin{Tm}
Given an $\mathbb{H}^{N\times N}$-valued function $K(p,q)$ on an
open set  $\Omega \subset \mathbb{H}$ let $\mathscr{H}(K)$ be
the associated reproducing kernel quaternionic Hilbert space.
Assume that for all $q\in \Omega$ the function
$p\mapsto K(p,q)$ is slice hyperholomorphic. Then the entries of
the elements of $\mathscr{H}(K)$ are also slice hyperholomorphic.
\end{Tm}
\begin{proof}
It is enough to consider the case of $\mathbb{H}$-valued function
because the matrix case works similarly. For any $f\in
\mathscr{H}(K)$, $p,q\in \Omega$ and $\varepsilon \in
\mathbb{R}\setminus \{0\}$ sufficiently small, we have
$$
\frac{1}{\varepsilon}(K(p,q+\varepsilon
)-K(p,q))=\frac{1}{\varepsilon}\overline{(K(q+\varepsilon
,p)-K(q,p))}.
$$
Let $(u+Iv, x+Iy)\in  \mathbb{C}_I\times \mathbb{C}_I$. We have that
$$
\frac{\partial K(p,q)}{\partial x}=\frac{\partial
\overline{K(q,p)}}{\partial u}.
$$
In an analogous way, we have:
$$
\frac{1}{\varepsilon}(K(p,q+I\varepsilon
)-K(p,q))=\frac{1}{\varepsilon}\overline{(K(q+I\varepsilon
,p)-K(q,p))},
$$
from which we deduce
$$
\frac{\partial K(p,q)}{\partial y}=\frac{\partial
\overline{K(q,p)}}{\partial v}.
$$
The two families
$$
\left\{\frac{1}{\varepsilon}(K(p,q+\varepsilon)-K(p,q))\right\}_{\varepsilon
\in \mathbb{R}\setminus \{0\}},  \qquad
\left\{\frac{1}{\varepsilon}(K(p,q+I\varepsilon)-K(p,q))\right\}_{\varepsilon
\in \mathbb{R}\setminus \{0\}},
$$
are uniformly bounded in the norm and therefore have weakly
convergent subsequences which converge to
$\displaystyle\frac{\partial K(p,q)}{\partial x}$ and
$\displaystyle\frac{\partial K(p,q)}{\partial y}$,  respectively.
Moreover we have
$$
\frac{1}{\varepsilon}(f(p+\varepsilon)-f(p))=\langle f(\cdot),
\frac{1}{\varepsilon}(K(\cdot,p+\varepsilon)-K(\cdot,p))\rangle_{\mathscr{H}(K)}
$$
and
$$
\frac{1}{\varepsilon}(f(p+I\varepsilon)-f(p))=\langle f(\cdot),
\frac{1}{\varepsilon}(K(\cdot,p+I\varepsilon)-K(\cdot,p))\rangle_{\mathscr{H}(K)}.
$$
Thus we can write
$$
\frac{\partial f}{\partial u}(p)=\langle f(\cdot), \frac{\partial
K(\cdot, p)}{\partial x} \rangle_{\mathscr{H}(K)},
$$
and
$$
\frac{\partial f}{\partial v}(p)=\langle f(\cdot), \frac{\partial
K(\cdot, p)}{\partial y} \rangle_{\mathscr{H}(K)}.
$$
To show that the function $f$ is slice hyperholomorphic, we
consider its restriction to any complex plane $\mathbb{C}_I$ and
we show that it is in the kernel of the corresponding
Cauchy-Riemann operator:
\[
\begin{split}
\frac{\partial f}{\partial u}+I\frac{\partial f}{\partial v}&=
\langle f, \frac{\partial K(\cdot, q)}{\partial x}
\rangle_{\mathscr{H}(K)}
+I \langle f(\cdot), \frac{\partial K(\cdot, q)}{\partial y} \rangle_{\mathscr{H}(K)}\\
  &= \langle f, \frac{\partial K(\cdot, q)}{\partial x}-\frac{\partial K(\cdot, q)}{\partial y}I\rangle_{\mathscr{H}(K)}\\
  &= \langle f, \overline{\frac{\partial K(q,\cdot)}{\partial u}+I\frac{\partial K(q,\cdot)}{\partial v}}\rangle_{\mathscr{H}(K)}=0\\
\end{split}
\]
since the kernel $K(q,p)$ is slice hyperholomorphic in the first
variable $q$.
\end{proof}


\section{Blaschke products}
\setcounter{equation}{0}
\label{secHardy}
The space $\mathbf{H}_2(\mathbb{B})$ was introduced in \cite{acs1} as
the space of power series $f(p)=\sum_{n=0}^\infty p^nf_n$, where
the coefficients $f_n\in\mathbb H$ and are such that
\begin{equation}
\label{normH2}
\|f\|_{\mathbf H_2(\mathbb B)}\stackrel{\rm
def.}{=} \sqrt{\sum_{n=0}^\infty |f_n|^2}<\infty.
\end{equation}
$\mathbf{H}_2(\mathbb{B})$ endowed with the inner
product
\[
[f,g]_2=\sum_{n=0}^\infty \overline{g_n}f_n,\quad{\rm where}\quad
g(p)=\sum_{n=0}^\infty p^ng_n
\]
is the right quaternionic reproducing kernel Hilbert space with
reproducing kernel
$$
k(p,q)=\sum_{n=0}^\infty p^n \bar q^n=(1-p\overline{q})^{-\star}.
$$
 The norm
\eqref{normH2} admits another expression.

\begin{Tm}
 The norm in $\mathbf{H}_2(\mathbb{B})$ can be written as
$$
\sup_{0<r<1,\ \ I\in\mathbb{S}}\Big[\frac{1}{2\pi}\int_0^{2\pi}
|f(re^{I\theta})|^2\, d\theta\Big]^{1/2}=
\sup_{0<r<1}\Big[\frac{1}{2\pi}\int_0^{2\pi}
|f(re^{I\theta})|^2\, d\theta\Big]^{1/2}.
$$
\end{Tm}
\begin{proof}
 When one writes the power series expansion for $f$ with center at $0$,
the equality is clear by the Parseval identity. Thus the norm can
be defined as in the classical complex case by computing the integral on a
chosen complex plane.
\end{proof}
Let us prove some results associated to the Blaschke factors
in the slice hyperholomorphic setting.
\begin{Dn}
Let $a\in\mathbb{H}$, $|a|<1$. The function
\begin{equation}
\label{eqBlaschke}
B_a(p)=(1-p\bar a)^{-\star}\star(a-p)\frac{\bar a}{|a|}
\end{equation}
is called  Blaschke factor at $a$.
 \end{Dn}
\begin{La}
\label{laba}
Let $a\in\mathbb B$. Then,
$B_a(p)$ is a slice
hyperholomorphic function in $\mathbb{B}$. Furthermore it holds that
\begin{equation}
\label{Baabar}
B_a(\overline{a})\overline{a}=\overline{a} B_a(\overline{a}).
\end{equation}
\end{La}

\begin{proof}
Indeed $B_a(p)$ is slice hyperholomorphic by its definition, moreover we have
\begin{equation}
\begin{split}
B_a(p)&=(\sum_{n=0}^\infty
p^n\overline{a}^n)\star(a-p)\frac{\overline{a}}{|a|}\\
&=\sum_{n=0}^\infty
(p^n\overline{a}^na-p^{n+1}\overline{a}^n)\frac{\overline{a}}{|a|}\\
&=|a|+\sum_{n=0}^\infty
p^{n+1}\overline{a}^{n+1}(|a|-\frac{1}{|a|}).
\end{split}
\label{Ba}
\end{equation}
Finally, \eqref{Baabar} is a direct consequence of the last equality.
\end{proof}

\begin{Rk}\label{Bpointwise}{\rm
Set $\lambda(p)=1-p\bar a$. Then
$$
(1-p\bar
a)^{-\star}=(\lambda^c(p)\star \lambda(p))^{-1}\lambda^c(p).
$$
Applying formula (\ref{productfg}) to the products $\lambda^c(p)\star\lambda(p)$ and $\lambda^c(p)\star (a-p)$,
we can rewrite (\ref{eqBlaschke}) as
\begin{equation}\label{lineartransf}
\begin{split}
B_a(p)&=(\lambda^c(p)\star \lambda(p))^{-1}\lambda^c(p)\star (a-p)\frac{\bar a}{|a|}=(\lambda^c(p)\lambda (\tilde p))^{-1}\lambda^c(p)(a-\tilde p)\frac{\bar a}{|a|}\\
&=\lambda (\tilde p)^{-1}(a-\tilde p)\frac{\bar a}{|a|}=(1-\tilde{p}\bar a)^{-1}(a-\tilde p)\frac{\bar a}{|a|},\\
\end{split}
\end{equation}
where $\tilde p=\lambda^c(p)^{-1} p \lambda^c(p)$.  Formula (\ref{lineartransf}) represents the Blaschke factor $B_a(p)$ in terms of pointwise multiplication only.}
\end{Rk}

\begin{Tm}\label{elrhfn}
Let $a\in\mathbb{H}$, $|a|<1$. The Blaschke factor $B_a(q)$ has
the following properties:
\begin{enumerate}
\item
it takes the unit ball $\mathbb{B}$ to itself;
\item it takes the boundary of the unit ball to itself;
\item it has a unique zero for $p=a$.
\end{enumerate}
\end{Tm}
\begin{proof} By Remark \ref{Bpointwise} we write
 $B_a(p)=(1-\tilde p\bar a)^{-1}(a-\tilde p)\dfrac{\bar a}{|a|}$.
 Let us show
that $|p|=|\tilde p|<1$ implies  $|B_a(p)|^2<1$. The latter inequality is
equivalent to
$$
|a-\tilde  p|^2<|1-\tilde p\bar a|^2
$$
which is also equivalent to
\begin{equation}\label{dffiuygah}
|a|^2+|p|^2<1+|a|^2|p|^2.
\end{equation}
The inequality (\ref{dffiuygah}) can be written as $(|p|^2-1)(1-|a|^2)<0$ and it holds when $|p|<1$.
When $|p|=1$ we set $p=e^{I\theta}$, so that $\tilde p=e^{I'\theta}$ by the proof of Corollary \ref{rotation}; we have
$$
|B_a(e^{I\theta})|=|1-e^{I'\theta}\bar
a|^{-1}|a-e^{I'\theta}|\frac{|\bar a|}{|a|}=|e^{-I'\theta}-\bar
a|^{-1}|a-e^{I'\theta}|=1.
$$
Finally, from (\ref{lineartransf}) it follows that  $B_a(p)$ has only one zero that comes from the factor $a-\tilde p$. Moreover
$B_a(a)=(1-\tilde{a}\bar a)^{-1}(a-\tilde a)\frac{\bar a}{|a|}$
where $\tilde a= (1-a^2)^{-1}a (1-a^2)=a$ and thus $B_a(a)=0$.
\end{proof}

\begin{Tm}\label{converge}
Let $\{a_j\}\subset \mathbb{B}$, $j=1,2,\ldots$ be a sequence
of nonzero quaternions such that $[a_i]\not=[a_j]$ if $i\not=j$ and assume that  $\sum_{j\geq 1} (1-|a_j|)<
\infty$. Then the function
\begin{equation}\label{B_product}
B(p):=\Pi^\star_{j\geq 1}(1-p\bar
a_j)^{-\star}\star(a_j-p)\frac{\bar a_j}{|a_j|},
\end{equation}
where $\Pi^\star$ denotes the $\star$-product, converges
uniformly on the compact subsets of $\mathbb{B}$.
\end{Tm}
\begin{proof}
Let $\alpha_j(p):=
B_{a_j}(p)-1$.
Using Remark \ref{Bpointwise} we have the chain of equalities:
\[
\begin{split}
\alpha_j(p)=& B_{a_j}(p)-1=(1-\tilde p \bar a_j)^{-1}(a_j-\tilde p)\frac{\bar a_j}{|a_j|}-1\\
=& (1-\tilde p \bar a_j)^{-1}\left[(a_j-\tilde p)\frac{\bar a_j}{|a_j|}-(1-\tilde p \bar a_j)\right]\\
=& (1-\tilde p \bar a_j)^{-1}\left[(|a_j|-1)\left(1+\tilde p\frac{\bar a_j}{|a_j|}\right)\right].
\end{split}
\]
Thus, if $|p|<1$ and recalling that $|\tilde p|=|p|$, we have
$$
|\alpha_j(p)|\leq 2 (1-|p|)^{-1}(1-|a_j|)
$$
and since $\sum_{j= 1}^\infty (1-|a_j|)<
\infty$ then $
\sum_{j=1}^\infty |\alpha_j(p)|
$
converges in $\mathbb{B}$ and the statement follows.
\end{proof}
\begin{Dn}
The function $B(p)$ defined in (\ref{B_product}) is called Blaschke product.
\end{Dn}
\begin{Rk}{\rm
In the complex case the sequence of complex numbers $\{a_j\}$ turns out to be the sequence of
zeroes of the Blaschke product. In the quaternionic case the
situation is different and we shall discuss it in the next results.
In order to illustrate the differences with the complex case, let us consider the simpler case in which we have a polynomial
$$
P(p)=(p-a_1)\star\ldots \star (p-a_n)
$$
and assume that $[a_i]\not=[a_j]$ for all $i,j=1,\ldots,n$. Then, it can be verified that
$p=a_1$ is a zero for the polynomial $P(p)$ while the other
zeroes belong to the spheres $[a_j]$ defined by $a_j$ for
$j=2,\ldots, n$. Note that, in the case in which  all the elements $a_j$
belong to a same sphere for all $j=1,...,n$, then the only zero
of the polynomial is $a_1$, see \cite[Lemma 2.2.1]{PhDPereira} and it has multiplicity $n$.
Moreover, whenever a polynomial and, more in general, a slice hyperholomorphic
function $f$ has two zeroes belonging to a same 2-sphere, then
all the elements of the sphere are zeroes for $f$.   Thus the zeroes of a slice hyperholomorphic
function are either isolated points or isolated
spheres, see \cite{MR2752913}.}
\end{Rk}

Assume that the slice hyperholomorphic function $f$ has zero set
$$
Z=\{a_1, a_2, \ldots\}\cup \{[c_1], [c_2], \ldots\}.
$$
Then it is possible to construct a suitable Blaschke product
having $Z_f$ as zero set. Let us begin with the case in which the
zeros are isolated points.
In the sequel, we will be in need of the following remark:
\begin{Rk}\label{commute}{\rm
Direct computations show the following equality of polynomials:
$$
(1-pa)\star (a-p)\frac{\bar a}{|a|}=\Big((a-p)\frac{\bar a}{|a|}\Big)\star (1-pa)=(a-p)\star (1-pa)\frac{\bar a}{|a|}.
$$
}
\end{Rk}

\begin{Pn}\label{isolated} Let $Z=\{a_1, a_2, \ldots\}$ be a sequence of elements in $\mathbb{B}$, $a_j\not=0$ for all $j=1,2,\ldots $ such that $[a_i]\not=[a_j]$ if $i\not=j$ and assume that $\sum_{j\geq 1}(1-|a_j|)<\infty$. Then there exists a Blaschke product $B(p)$ having zero set at $Z$.
\end{Pn}
\begin{proof}
Let us prove the statement by induction. By hypothesis the zero set
of the required Blaschke product consists of isolated points, all
of them belonging to different spheres. If $n=1$, we have already
proved that $B_1(p):=B_{a_1}(p)$ has $a_1$ as its unique zero.
Let us assume that the statement holds for $a_1,\ldots ,a_k$, and
so there exists a Blaschke product $B_k(p)$ vanishing at the
given points and  let us prove that we can construct a Blaschke
product vanishing at $a_1,\ldots ,a_k, a_{k+1}$. Observe that it is possible
to choose an element $a'_{k+1}$  belonging to the sphere
$[a_{k+1}]$ such that
$$
B_{k}(p)\star B_{a'_{k+1}}(p)
$$
has zeros $a_1,\ldots, a_{k+1}$. In fact, consider the product
$$
B_{k+1}(p):=B_k(p)\star (1-p\bar
a'_{k+1})^{-\star}\star(a'_{k+1}-p)\frac{\bar a_{k+1}'}{|a_{k+1}'|}
$$
and rewrite it using Remark \ref{commute} in the form
$$
B_{k+1}(p):=B_k(p)\star (a'_{k+1}-p)\star (1-p
a'_{k+1})(1-2{\rm Re}(\bar
a'_{k+1})p+|\bar
a'_{k+1}|^2p^2)^{-1}\frac{\bar a_{k+1}'}{|a_{k+1}'|}.
$$
We now observe that the zeros of $B_{k+1}(p)$ belonging to the ball $\mathbb{B}$ come from the zeros of the product
$$
\tilde{B}(p):=B_k(p)\star (a'_{k+1}-p).
$$
Observe that
 $$
\tilde{B}(a_{k+1})=  B_k(a_{k+1})(a'_{k+1}-B_k(a_{k+1})^{-1}a_{k+1}B_k(a_{k+1}))
 $$
 and in order that $a_{k+1}$ is a zero of $\tilde{B}$, and so of $B_{k+1}$, it is sufficient to choose
 $$
 a'_{k+1}=B_k(a_{k+1})^{-1}a_{k+1}B_k(a_{k+1}).
 $$
 The convergence of the Blaschke product follows as in Theorem \ref{converge}.
 \end{proof}

From now on, when we write $Z=\{(a,\mu)\}$ we mean that $Z$ consists of the point $a$ repeated $\mu$ times.
Let us now prove the analog of Theorem \ref{elrhfn} (3) in the case in which the point $a$ has multiplicity $\mu$.
 \begin{La}\label{multiple}
 Let $Z=\{(a,\mu)\}$ with $a\in \mathbb{B}$ and $a\not=0$. The Blaschke product
 \[
 \begin{split}
B(p):=& \Big((1-p\bar
a)^{-\star}\star(a-p)\frac{\bar a}{|a|}\Big)^{\star\mu}:=\\
=&\underbrace{\begin{pmatrix}(1-p\bar
a)^{-\star}\star(a-p)\frac{\bar a}{|a|}& \ldots & (1-p\bar
a)^{-\star}\star(a-p)\frac{\bar a}{|a|}\end{pmatrix}}_{\mu\,\,{\rm
times}}\\
\end{split}
\]
has $Z$ as zero set.
\end{La}
\begin{proof}
We have:
$$
(1-p\bar
a)^{-\star}\star(a-p)\frac{\bar a}{|a|}=(1-2{\rm Re}(a)p +p^2|a|^2)^{-1}(1-pa)\star(a-p)\frac{\bar a}{|a|}
$$
thus, using the fact that $1-2{\rm Re}(a)p +p^2|a|^2$ has real coefficients, we can write
\[
B(p)=(1-2p{\rm Re}(a) +p^2|a|^2)^{-\mu}\underbrace{\begin{pmatrix}(1-pa)\star(a-p)\frac{\bar a}{|a|} &\ldots &
(1-pa)\star(a-p)\frac{\bar a}{|a|}\end{pmatrix}}_{\mu\,\,{\rm
times}}
\]
and thanks to Remark \ref{commute} we obtain
$$
B(p)=(1-2p{\rm Re}(a) +p^2|a|^2)^{-\mu}\Big((a-p)\frac{\bar a}{|a|}\Big)^\mu\star (1-pa)^\mu.
$$
Thus $B(p)$ has, in $\mathbb{B}$, a unique zero at $p=a$ of multiplicity $\mu$. Note that the zero on the sphere $[1/a]$ which, as it can be proven,
coincides with $1/a$ has to be excluded since $B(p)$ is not defined there, moreover $1/a\not\in \mathbb{B}$.
\end{proof}
\begin{Pn}\label{sfericalzeros} Let $Z=\{(a_1,\mu_1), (a_2,\mu_2), \ldots\}$ be a sequence of points $a_j\in \mathbb{B}$ with  respective multiplicities $\mu_j\geq 1$, $a_j\not=0$ for $j=1,2, \ldots$. Let $a_j$ be such that $[a_i]\not=[a_j]$ if $i\not=j$ and $\sum_{j\geq 1} \mu_j (1-|a_j|)<
\infty$. Then there exists a Blaschke product of the form
$$
B(p)=\prod_{j\geq 1}^\star (B_{a'_j}(p))^{\star \mu_j},
$$
 having zero set at $Z$,
where $a'_1=a_1$ and $a'_j\in [a_j]$ are suitably chosen elements, $j=2,3,\ldots$.
\end{Pn}
\begin{proof}
We prove the assertion by induction on the number of distinct zeros. If there is just one zero $a_1$ with multiplicity $\mu_1$, then
the statement follows by Lemma \ref{multiple}.
Let us assume that the statement holds in the case we have $k$ different zeros $a_i$ with respective multiplicities $\mu_i$ and let us prove that it holds for $k+1$ different zeros.
Let $B_k(p)$ be the Blaschke product having zeros at $Z=\{(a_1,\mu_1), \ldots,(a_k,\mu_k)\}$ and let us consider
$$
B_{k+1}(p):=B_k(p)\star (B_{a'_{k+1}}(p))^{\star \mu_k}
$$
where $a'_{k+1}$ is chosen such that $B_k(p)\star B_{a'_{k+1}}(p)$ has a zero at $p=a_{k+1}$. Then all the other zeros of $B_{k+1}$ must belong to the sphere $[a_{k+1}]$. Moreover they must coincide with $a_{k+1}$ otherwise the Blaschke product $(B_{a'_{k+1}}(p))^{\star \mu_k}
$ vanishes at two different points on a same sphere, and thus it vanishes on the whole sphere. In particular, any two conjugate elements on the sphere are zeros of the product and so we would have:
 \[
 \begin{split}
 &B_a(p)\star B_{\bar a}(p)=(1-p\bar a)^{-\star}\star(a-p)\frac{\bar a}{|a|}\star (1-pa)^{-\star}\star(\bar a-p)\frac{a}{|a|}\\
 &=(1-2{\rm Re}(a)p+p^2|a|^2)^{-1}(|a|^2-2{\rm Re}(a)p+p^2).
 \end{split}
 \]
However, it is immediate that the product $(B_{a'_{k+1}}(p))^{\star \mu_k}$ does not contain  factors of the above form, thus all its zeros coincide with $a_{k+1}$ as stated.
 The convergence of the Blaschke product follows as in Theorem \ref{converge}.
\end{proof}
 If a Blaschke product of two factors has an entire sphere of zeros then, as discussed in the proof of the previous theorem, it has a specific form and we are led to the following definition:
 \begin{Dn}
Let $a\in\mathbb{H}$, $|a|<1$. The function
\begin{equation}
\label{blas_sph}
B_{[a]}(p)=(1-2{\rm Re}(a)p+p^2|a|^2)^{-1}(|a|^2-2{\rm Re}(a)p+p^2)
\end{equation}
is called  Blaschke factor at the sphere $[a]$.
 \end{Dn}
 \begin{Rk}{\rm Note that the definition of $B_{[a]}(p)$ does not depend on the choice of the point $a$ that identifies the 2-sphere. Indeed all the elements in the sphere $[a]$ have the same real part and module.
 It is easy to verify that Blaschke factor $B_{[a]}(p)$ vanishes on the sphere $[a]$.}
 \end{Rk}
The  following result is immediate:
\begin{Pn} A Blaschke product having zeros at the set
of spheres
$$
Z=\{([c_1],\nu_1), ([c_2],\nu_2), \ldots\}
$$
where $c_j\in\mathbb{B}$, the sphere $[c_j]$ is a zero of multiplicity $\nu_j$, $j=1,2,\ldots$
and $\sum_{j\geq 1} \nu_j(1-|c_j|)<\infty$ is given by
\[
\prod_{j\geq 1}(B_{[c_j]}(p))^{\nu_j}.
\]
\end{Pn}
\begin{proof}
All the factors
 $B_{[c_j]}(p)$ have real coefficients and thus belong to the class $\mathcal{N}$ (see Section \ref{positivity}), so
 we can use the pointwise product. The fact that the zeros are the given spheres follows by taking the zeros of each factor.
 The convergence of the infinite product follows as in Theorem \ref{converge}.
\end{proof}
\begin{Tm}
A Blaschke product having zeros at the set
 $$
 Z=\{(a_1,\mu_1), (a_2,\mu_2), \ldots, ([c_1],\nu_1), ([c_2],\nu_2), \ldots \}
 $$
 where $a_j\in \mathbb{B}$, $a_j$ have respective multiplicities $\mu_j\geq 1$, $a_j\not=0$ for $j=1,2,\ldots $, $[a_i]\not=[a_j]$ if $i\not=j$, $c_i\in \mathbb{B}$, the spheres $[c_j]$ have respective multiplicities $\nu_j\geq 1$,
 $j=1,2,\ldots$, $[c_i]\not=[c_j]$ if $i\not=j$
and
$$
\sum_{i,j\geq 1} \Big(\mu_i (1-|a_i|)+
\nu_j(1-|c_j|)\Big)<\infty
$$
is given by
\[
\prod_{i\geq 1} (B_{[c_i]}(p))^{\nu_i}\prod_{j\geq 1}^\star (B_{a'_j}(p))^{\star \mu_j},
\]
where $a'_1=a_1$ and $a'_j\in [a_j]$ are suitably chosen elements, $j=2,3,\ldots$.
\end{Tm}
\begin{proof}
It follows from Propositions \ref{isolated} and \ref{sfericalzeros}.
\end{proof}
\begin{Tm}
\label{tmiso} Let $B_a$ be a Blaschke factor. The operator
\[
M_a\,:\, \,f\mapsto B_a\star f
\]
is an isometry from $\mathbf H_2(\mathbb B)$ onto itself.
\end{Tm}

\begin{proof} We first consider $f(p)=p^uh$ and $g(p)=p^vk$ where
$u,v\in\mathbb N_0$ and $h,k\in\mathbb H$, and show that
\begin{equation}
\label{eq:iso}
[B_a\star f,B_a\star g]_2=\delta_{uv}\overline{k}h.
\end{equation}

Using calculation (\ref{Ba}), and with $f$ and $g$ as above, we have
\[
(B_a\star f)(p)=p^uh|a|+\sum_{n=0}^\infty
p^{n+1+u}\overline{a}^{n+1}(|a|-\frac{1}{|a|})h
\]
and
\begin{equation}
\label{bag}
(B_a\star g)(p)=p^vk|a|+\sum_{n=0}^\infty
p^{n+1+v}\overline{a}^{n+1}(|a|-\frac{1}{|a|})k.
\end{equation}
If $u=v$ we have
\[
[B_a\star f,B_a\star g
]_2=\overline{k}h\left(|a|^2+\sum_{n=0}^\infty
|a|^{2n+2}(|a|-\frac{1}{|a|})^2\right)=\overline{k}h=[f,g]_2.
\]
To compute $[f,g]_2$ we assume that $u<v$.  Then, in view of
\eqref{bag} we have
\[
[p^uh|a|,B_a\star g]_2=0.
\]
So
\[
\begin{split}
[B_a\star f,B_a\star g]_2&=[\sum_{n=0}^\infty
p^{n+1+u}\overline{a}^{n+1}\left(|a|-\frac{1}{|a|}\right)h,p^v|a|k]_2+\\
&\hspace{5mm}+ [\sum_{n=0}^\infty
p^{n+1+u}\overline{a}^{n+1}\left(|a|-\frac{1}{|a|}\right)h,
\sum_{m=0}^\infty
p^{m+1+v}\overline{a}^{m+1}\left(|a|-\frac{1}{|a|}\right)k]_2\\
&=|a|\overline{k}\overline{a}^{v-u}\left(|a|-\frac{1}{|a|}\right)h+\\
&\hspace{5mm}+[\sum_{m=0}^\infty
p^{m+1+v}\overline{a}^{m+1+v-u}\left(|a|-\frac{1}{|a|}\right)h,
\sum_{m=0}^\infty
p^{m+1+v}\overline{a}^{m+1}\left(|a|-\frac{1}{|a|}\right)k]_2\\
&=|a|\overline{k}\overline{a}^{v-u}\left(|a|-\frac{1}{|a|}\right)h
+\overline{k}\left(|a|-\frac{1}{|a|}\right)^2\overline{a}^{v-u}
\frac{|a|^2}{1-|a|^2}h\\
&=0\\
&=[f,g]_2.
\end{split}
\]
The case $v<u$ is considered by symmetry of the inner product.
Hence, \eqref{eq:iso} holds for polynomials. By continuity, and a corollary of the Runge theorem, see \cite{2763766}, it
holds for all $f\in\mathbf H_2(\mathbb B)$.

\end{proof}

We mention that similar computations hold in the case of
bicomplex numbers. See \cite{alss2}.

\section{Interpolation in the Hardy space}
\setcounter{equation}{0}
In this section we consider the following problem:
\begin{Pb}
\label{pb1}
Given $N$ points $a_1,\ldots ,a_N\in\mathbb{B}$, and $M$ spheres $[c_1], \ldots, [c_M]$ in $\mathbb B$ such that the spheres $[a_1],\ldots [a_N],
[c_1], \ldots, [c_M]$  are pairwise non-intersecting,
find all $f\in\mathbf H_2(\mathbb B)$ such that
\begin{equation}
\label{eq:cond1}
f(a_i)=0,\quad i=1,\ldots, N,
\end{equation}
and
\begin{equation}
f([c_j])=0,\quad j=1,\ldots, M.
\label{eq:cond2}
\end{equation}
\end{Pb}
\begin{Tm}
\label{tminter}
There is a Blaschke product $B$ such that the
solutions of Problem \ref{pb1} are the functions $f=B\star g$,
when $g$ runs through $\mathbf H_2(\mathbb B)$.
\end{Tm}

We give two proofs of this theorem, the first one iterative, using
formula \eqref{productfg}, and the second one global.

\begin{proof} (Iterative proof). We proceed in three steps.
As a preliminary computation
we consider in the first step the case $N=1$ and $M=0$. The problem itself will
be solved by considering the interpolation at the spheres first.\\

STEP 1: {\sl We solve the problem for $M=0$ and $N=1$.}\\

Let $B_{a_1}$ be the Blaschke factor \eqref{eqBlaschke}
at $a_1$. By \eqref{productfg}, we
have $(B_{a_1}\star f)(a_1)=0$ for all $f\in\mathbf H_2(\mathbb
B)$. Furthermore, by Theorem \ref{tmiso} we have that
$\|B_{a_1}\star f\|_2=\|f\|_2$. Thus for $N=1$ the set $\mathcal
M$ of solutions to Problem \ref{pb1} contains $B_{a_1}\star
\mathbf H_2(\mathbb B)$. We now prove that $\mathcal M\subseteq
B_{a_1}\star\mathbf H_2(\mathbb B)$. Let $f\in\mathcal M$. Then, by the
reproducing kernel property, $f$ is orthogonal to
$(1-p\overline{a_1})^{-\star}$. The range ${\rm
ran}~\sqrt{I-M_{a_1}M_{a_1}^*}$ is equal to the span of
$(1-p\overline{a_1})^{-\star}$ (see \cite{acs1}). In view of
Theorem \ref{tmiso} we have
$\sqrt{I-M_{a_1}M_{a_1}^*}=I-M_{a_1}M_{a_1}^*$ and thus:
\[
\mathbf H_2(\mathbb B)=(I-M_{a_1}M_{a_1}^*)\mathbf H_2(\mathbb
B)\oplus (M_{a_1}M_{a_1}^*)\mathbf H_2(\mathbb B).
\]
Therefore $f\in (M_{a_1}M_{a_1}^*)\mathbf H_2(\mathbb B)$. Hence
$\mathcal M=B_{a_1}\star \mathbf H_2(\mathbb B)$.\\

With this preliminary computation at hand, we solve the
interpolation problem by first considering the interpolation at
the spheres
$[c_1], \ldots, [c_M]$.\\

STEP 2: {\sl  Consider the sphere $[c_j]$ and let $B_{[c_j]}$ be the
corresponding Blaschke factor given by
\eqref{blas_sph}, $j=1,2, \ldots, M$. An element $f\in\mathbf
H_2(\mathbb B)$ vanishes on the spheres $[c_1],\ldots , [c_M]$ if and only it
can be written as
\begin{equation}
\label{interbl}
f=B_{[c_1]}B_{[c_2]}\cdots B_{[c_M]}g,
\end{equation}
where $g\in\mathbf H_2(\mathbb B)$.}\\

Note that in \eqref{interbl} we have pointwise products since the
Blaschke factors on spheres have real coefficients.   By \cite[Corollary 4.3.7, p. 123]{MR2752913}, $f$
vanishes on the whole sphere $[c_1]$ if and only if
$f(c_1)=f(\overline{c_1})=0$. By STEP 1, the first condition
means that $f=B_{c_1}\star g$ for some $g\in\mathbf H_2(\mathbb B)$. By
\eqref{productfg}, the second condition is equivalent to:
\begin{equation}
\label{c1bar}
B_{c_1}(\overline{c_1})
g((B_{c_1}(\overline{c_1}))^{-1}\overline{c_1}B_{c_1}(\overline{c_1}))=0.
\end{equation}
Since $B_{c_1}(\overline{c_1})\not=0$, and taking into account
\eqref{Baabar}, we see that \eqref{c1bar} is equivalent to
$g(\overline{c_1})=0$. Thus, once more using STEP 1, we have
$g(p)=B_{\overline{c_1}}\star h$ for some $h \in\mathbf
H_2(\mathbb B)$. Therefore
\[
f=B_{c_1}\star B_{\overline{c_1}}\star h=B_{[c_1]}h.
\]

This argument can be iterated for the spheres $[c_2],\ldots,
[c_M]$ since $B_{[c_2]}(c_1)\not=0$ (this last inequality in
turn following from the fact that the spheres do not intersect).\\

We now turn to the conditions \eqref{eq:cond1}. The function $f$
is of the form \eqref{interbl}, and thus the condition $f(a_1)=0$
becomes
\[
\left(B_{[c_1]}B_{[c_2]}\cdots B_{[c_M]}\right)(a_1)g(a_1)=0,
\]
and so, by STEP 1, $g=B_{a_1}\star g_1$ for some $g_1\in\mathbf
H_2(\mathbb B)$. Let now $f\in\mathbf H_2(\mathbb B)$ satisfying
moreover $f(a_2)=0$. By the previous argument, $f$ is of the form
\[
\left(B_{[c_1]}B_{[c_2]}\cdots B_{[c_M]}\right)B_{a_1}\star g_1
\]
for some $g_1\in\mathbf H_2(\mathbb B)$. The condition $f(a_2)=0$
and formula \eqref{productfg} gives
\[
g(a_2^\prime)=0, \quad{\rm where}\quad a_2^\prime=X^{-1}a_2X,
\]
with
\[
X=\left(B_{[c_1]}B_{[c_2]}\cdots B_{[c_M]}B_{a_1}\right)(a_2).
\]
Hence $f$ is a solution if and only if $g_2=B_{a_2^\prime}\star
g_2$ for some $g_1\in\mathbf H_2(\mathbb B)$. The argument can be
iterated and we obtain the set of all functions $f\in\mathbf
H_2(\mathbb B)$ which vanish at the points $a_1,\ldots, a_N$.
\end{proof}

We now turn to the global proof of Theorem \ref{tminter}.

\begin{proof} (Global proof).
We proceed in a number of steps and first define $\mathcal M$ to
be the span of the functions
\[
(1-pa_1)^{-\star},\ldots,
(1-pa_N)^{-\star},(1-pc_1)^{-\star},(1-p\overline{c_1})^{-\star},\ldots,
(1-pc_M)^{-\star},(1-p\overline{c_M})^{-\star}.
\]
Define
\[
A={\rm diag}~(a_1,\ldots,
a_N,c_1,\overline{c_1},\ldots,c_M,\overline{c_M}),
\]
and
\[
c=\underbrace{\begin{pmatrix}1&1&\cdots&1&1\end{pmatrix}}_{(N+2M)\,\,{\rm
times}}.
\]
Finally let $\mathbf P$ denote the Gram matrix of $\mathcal M$ in
the $\mathbf H_2(\mathbb B)$ inner product. The claim of the first
step is a direct consequence of the reproducing kernel property in
$\mathbf H_2(\mathbb B)$.\\

STEP 1: {\sl $f\in\mathbf H_2(\mathbb B)$ is a solution of the
interpolation problem \ref{pb1} if and only if it is orthogonal
to $\mathcal M$.}\\

STEP 2: {\sl The matrix $\mathbf P\in\mathbb H^{(N+2M)\times
(N+2M)}$ is strictly positive and satisfies the matrix equation
(called a Stein equation)}
\begin{equation}
\label{stein}
\mathbf P-A^*\mathbf P A=c^*c.
\end{equation}

This step is also a consequence of the reproducing kernel
property since, for $a,b\in\mathbb B$ it holds that:
\[
[(1-p\overline{b})^{-\star},(1-p\overline{a})^{-\star}]_2=\sum_{n=0}^\infty
a^n\overline{b}^n,
\]
and so
\[
[(1-p\overline{b})^{-\star},(1-p\overline{a})^{-\star}]_2-a[
(1-p\overline{b})^{-\star},(1-p\overline{a})^{-\star}]_2\overline{b}=1.
\]
\mbox{}\\

STEP 3: {\sl There exists a vector $b\in\mathbb H^{N+2M}$ and
$d\in\mathbb H$ such that}
\begin{equation}
\begin{pmatrix}A&b\\ c&d\end{pmatrix}\begin{pmatrix} \mathbf P^{-1}&0\\
0&1\end{pmatrix}\begin{pmatrix}A&b\\ c&d\end{pmatrix}^*=
\begin{pmatrix} \mathbf P^{-1}&0\\
0&1\end{pmatrix}.
\end{equation}

Since $\mathbf P>0$, it has a strictly positive square root. For
this last fact, see for instance \cite[Proposition 3.1.3, p. 440]{as3}, and see
the references \cite[Corollary 6.2, p. 41]{MR97h:15020}
\cite{MR13:312d, MR12:153i} for the structure of
normal quaternionic matrices. Thus, the Stein equation
\eqref{stein} can be rewritten as $T^*T=I_{N+2M}$, where
\[
T=\begin{pmatrix} \mathbf P^{1/2}A\mathbf P^{-1/2}\\
c\mathbf P^{-1/2}
\end{pmatrix}.
\]
Thus the columns of $T$ form an $N+2M$ orthogonal sets of vectors
in $\mathbb H^{N+2M+1}$, and we can complete it to an orthonormal
basis with a vector $h\in\mathbb H^{N+2M+1}$ such that
\[
\begin{pmatrix}T&h\end{pmatrix}^*\begin{pmatrix}T&h\end{pmatrix}=I_{N+2M+1}.
\]
The claim follows from the fact that $TT^*=I_{N+2M+1}$ with
\begin{equation}
\label{eq:bd}
\begin{pmatrix}b\\d\end{pmatrix}=\begin{pmatrix}\mathbf
P^{-1/2}&0\\0&1\end{pmatrix}h.
\end{equation}

We now introduce

\[
\begin{split}
B(p)=d+ pc \star (I-pA)^{-\star} b
=d+\sum_{n=0}^\infty p^{n+1}cA^nb.
\end{split}
\]

Step 4 below is a particular case of Proposition \ref{ag} below, and its
proof will be omitted.\\

STEP 4: {\sl The function $B$ satisfies
\begin{equation}
\label{BSchurmult}
C(I_{N+2M}-pA)^{-\star}\mathbf P^{-1}(C(I_{N+2M}-qA)^{-\star})^*
=(1-B(p)\overline{B(q)})\star (1-p\overline{q})^{-\star}.
\end{equation}
}

Since $ \mathbf{P}^{-1}>0$ it follows from \eqref{BSchurmult} that
$B$ is a Schur multiplier, and in particular, ${\mathcal M}={\rm
ran}~\sqrt{I-M_BM_B^*}$, where $M_B$ denotes the operator of
slice multiplication by $B$ on the left (the
square root exists because the operator has finite rank; more
generally, any positive operator in a quaternionic Hilbert space
has a positive square root. We will not need this general fact
here). Since $\mathcal M$ is finite dimensional, we have more
precisely
\begin{equation}
{\mathcal M}={\rm ran}~\sqrt{I-M_BM_B^*}={\rm ran}~(I-M_BM_B^*).
\label{eq:range}
\end{equation}

STEP 5: {\sl The function $B\star g$ satisfies the interpolation conditions
\eqref{eq:cond1}-\eqref{eq:cond2} for every $g\in\mathbf H_2(\mathbb B)$.}\\

We first
prove that $(B\star g)(a_1)=0$. The proof that $B\star g$ vanishes at the
points $a_2,\ldots ,a_N$ and $c_1,\overline{c_1},\ldots,c_M,\overline{c_M}$
is the same.
\[
\begin{split}
B(a_1)&=d+\sum_{n=0}^\infty a_1^{n+1}\begin{pmatrix}\overline{a_1}^n&
\overline{a_2}^n&\cdots&\overline{a_N}^n\end{pmatrix}b\\
&=d+a_1\begin{pmatrix}1&0&\cdots &0\end{pmatrix}\mathbf P b\\
&=\begin{pmatrix}a_1\begin{pmatrix}1&0&\cdots &0\end{pmatrix}\mathbf P &1
\end{pmatrix}\begin{pmatrix}b\\d\end{pmatrix}\\
&=
\begin{pmatrix}a_1\begin{pmatrix}1&0&\cdots &0\end{pmatrix}\mathbf P^{1/2} &1
\end{pmatrix}h\qquad (\mbox{\rm where we have used \eqref{eq:bd}})\\
&=0
\end{split}
\]
by definition of $h$ since
\[
\begin{pmatrix}
\mathbf P^{1/2}
\begin{pmatrix}1\\ 0\\ \vdots \\0
\end{pmatrix}
\overline{a_1}\\1
\end{pmatrix}=T\left(\mathbf P^{1/2}
\begin{pmatrix}1\\ 0\\ \vdots \\0
\end{pmatrix}
\right)\in{\rm ran}~T.
\]
Let now $g(p)=p^mk$ with $m\in\mathbb N$ and $k\in\mathbb H$. Then,
\[
(B\star g)(p)=p^mB(p)k
\]
and so the interpolation conditions \eqref{eq:cond1}-\eqref{eq:cond2} are met.
The result is thus true for all slice hyperholomorphic polynomials in $p$, and
hence, in view of the preceding step, for every element $g$ in $\mathbf H_2(
\mathbb B)$ since convergence in norm implies pointwise convergence in a
reproducing kernel Hilbert space.\\

STEP 6: {\sl The set of solutions is given by $f=B\star g$, when
$g$ runs through $\mathbf H_2(\mathbb B)$.}\\

By definition of $\mathcal M$ and using \eqref{eq:range} and the reproducing
kernel property we see that ${\rm ran}~ (I-M_BM_B^*)$ and ${\rm ran}~ M_BM_B^*$
are orthogonal in $\mathbf H_2(\mathbb B)$. Since the sum of these two ranges
is the whole of $\mathbf H_2(\mathbb B)$, we deduce that
$\mathcal M^\perp={\rm ran}~M_BM_B^*$ and this ends the proof since
$
{\rm ran}~(M_BM_B^*)={\rm ran} ~M_B.
$
\end{proof}

\section{Quaternionic Pontryagin spaces}
\setcounter{equation}{0}
Quaternionic Pontryagin spaces have been studied in \cite{as3}.
In this section we review the main definitions, and prove in the
setting of quaternionic spaces, an important result due to
Shmulyan in the complex setting; see \cite{s} and \cite[Theorem
1.4.2, p. 29]{adrs}. Consider a right vector space $\mathscr P$ on
the quaternions, endowed with a $\mathbb H$-valued Hermitian form
$[\,\cdot\,,\, \cdot\,]$, meaning that
\[
[va,wb]=\overline{b}[v,w]a,\quad\forall a,b\in\mathbb H\quad{\rm
and}\quad\forall v,w\in \mathscr P.
\]
$\mathscr P$ is called a (right, quaternionic) Pontryagin space
if it admits a decomposition
\begin{equation}
\label{eq:sum}
\mathscr P=\mathscr P_++\mathscr P_-,
\end{equation}
into a sum of two vector subspaces $\mathscr P_+$ and $\mathscr
P_-$  with the following properties:\\
$(1)$ $(\mathscr P_+,
[\cdot,\cdot])$ is a (right, quaternionic)
Hilbert space.\\
$(2)$ $(\mathscr P_-,-[\cdot, \cdot])$ is a finite dimensional
(right, quaternionic)
Hilbert space.\\
$(3)$ The sum \eqref{eq:sum} is direct and orthogonal:
$\mathscr P_+\cap\mathscr P_-=\left\{0\right\}$ and
\[
[v_+,v_-]=0,\quad\forall v_+\in\mathscr P_+\quad{\rm
and}\quad\forall v_-\in\mathscr P_-.
\]
The space $\mathscr P$ endowed with the form
\[
\langle v,w\rangle=[v_+,v_-]-[w_+,w_-],\quad v=v_++v_-,\,\,
w=w_++w_-,
\]
is a (right quaternionic) Hilbert space. The decomposition
\eqref{eq:sum} is called a fundamental decomposition. It is not
unique (except for the case where one of the components reduces to
$\left\{0\right\}$), but all the corresponding Hilbert space
topologies are equivalent; see \cite[Theorem 12.3, p. 467]{as3}.
The number $\kappa={\rm dim}~\mathscr P_-$ is called the index of
the Pontryagin space $\mathscr P$. It is the same for all the decompositions;
see \cite[Proposition 12.6, p. 469]{as3}. The reader should be
aware that in some sources on the complex valued case, and in
particular in \cite{bognar, ikl}, the convention is the opposite,
and it is the space $\mathscr P_+$
which is required to be finite dimensional.\\

An important  example of finite dimensional Pontryagin space is:

\begin{Ex}
Let $J\in\mathbb H^{N\times N}$ be a signature matrix. The
space $\mathbb H^N$ endowed with the Hermitian form
\[
[v,w]_J=w^*Jv.
\]
is a right quaternionic Pontryagin space, which we will denote by
${\mathbb H}^{N}_{J}$.
\end{Ex}

Before turning to Shmulyan's theorem we recall the following
definition: Given two right quaternionic Pontryagin spaces
$(\mathscr P_1,[\cdot,\cdot]_1)$ and $(\mathscr
P_2,[\cdot,\cdot]_2)$ a linear relation between $\mathscr P_1$
and $\mathscr P_2$ is a right linear subspace, say $R$,  of the
product $\mathscr P_1\times\mathscr P_2$. The domain of the
relation is the set of elements $v_1\in\mathscr P_1$ such that
there exists a (not necessarily unique) $v_2\in\mathscr P_2$ such
that $(v_1,v_2)\in R$. The relation is called contractive if
\[
[v_1,v_1]_1\le [v_2,v_2]_2,\quad\forall (v_1,v_2)\in R.
\]
The graph of an operator is a relation. A relation will be the
graph of an operator if and only it has no elements of the form
$(0,v_2)$ with $v_2\not=0$.

\begin{Tm}
\label{schmu}
A densely defined contractive relation between
quaternionic Pontryagin spaces of same index extends to the graph
of a contraction from $\mathscr P_1$ into $\mathscr P_2$.
\end{Tm}

\begin{proof} We follow the strategy of \cite[p. 29-30]{adrs} we
divide the proof into a number of steps. We recall that a
 strictly negative subspace is a linear subspace
$V$ such that $[v,v]<0$ for every non zero element of $V$.\\

STEP 1: {\sl The domain of the relation contains a maximum
negative subspace.}\\

Indeed, every dense linear subspace of a right quaternionic
Pontryagin space of index $\kappa>0$ contains a $\kappa$
dimensional strictly negative subspace. See \cite[Theorem 12.8 p.
470]{as3}. We denote by $\mathscr V_-$ such a subspace of the
domain of $R$.\\

STEP 2: {\sl The relation $R$ restricted to $\mathscr V_-$ has a
zero kernel, and the image of $\mathscr V_-$ is a strictly
negative subspace of $\mathscr P_2$ of dimension $\kappa$.}\\

Let $(v_1,v_2)\in R$ with $v_1\in\mathscr V_-$. Since $R$ is
contractive  we have
\[
[v_2,v_2]_2\le [v_1,v_1]_1\le0,
\]
and the second inequality is strict when $v_1\not=0$. Thus, the
image of $\mathscr V_-$ is a strictly negative subspace of
$\mathscr P_2$. Next, let $(v,w)\in R$ and $(\widetilde{v},{w})$
with $v,\widetilde{v}\in\mathscr V_-$ and $w\in\mathscr P_2$.
Then, $(v-\widetilde{v},0)\in R$. Since $R$ is contractive we have
\[
[0,0]_2\le[v-\widetilde{v},v-\widetilde{v}]_1
\]
This forces $v=\widetilde{v}$ since $\mathscr V_-$ is strictly
negative, and proves the second step.\\

STEP 3: {\sl $R$ is the graph of a densely defined contraction.}\\

We choose $\mathscr V_-$ as in the first two steps, and take
$v_1,\ldots, v_\kappa$ a basis of $\mathscr V_-$. Then, there
exists uniquely defined vectors $w_1\,\ldots, w_\kappa\in\mathscr
P_2$ such that $(v_i,w_i)\in R$ for $i=1,\ldots, \kappa$. Set
$\mathscr W_-$ to be the linear span of $w_1,\ldots, w_\kappa$. By
Step 2 and since the spaces $\mathscr P_1$ and $\mathscr P_2$
have the same negative index
\[
{\rm dim}~\mathscr V_-={\rm dim}~\mathscr W_-={\rm ind}_-\mathscr
P_1={\rm ind}_-\mathscr P_2,
\]
and there exists fundamental decompositions
\[
\mathscr P_1=\mathscr V_-+\mathscr V_+\quad{\rm and}\quad
\mathscr P_2=\mathscr W_-+\mathscr W_+,
\]
where $(\mathscr V_+,[\cdot,\cdot]_1)$ and $(\mathscr
W_+,[\cdot,\cdot]_2)$ are right quaternionic Hilbert spaces. Let
now $(0,w)\in R$. We need to show that $w=0$. Still following
\cite[p. 30]{adrs} we write $w=w_-+w_+$ where $w_-\in\mathscr W_-$
and $w_+\in\mathscr W_+$. Let $w_-=\sum_{n=1}^\kappa w_jq_j$
where the $q_j\in\mathbb H$, and set $v_-=\sum_{n=1}^\kappa
v_jq_j$. Then, $(v_-,w_-)\in R$ and
\[
(0,w)=(v_-,w_-)+(-v_-,w_+).
\]
It follows that $(-v_-,w_+)\in R$. Since $R$ is contractive, we
have
\[
[w_+,w_+]_2\le [v_-,v_-]_1,
\]
and so $[w_+,w_+]_2\le 0$. Thus $w_+=0$. It follows that
$(0,w_-)\in R$ and so $w_-=0$, because $R$ is one-to-one on
$\mathscr V_-$, as follows from STEP 2.\\

STEP 4: {\sl $R$ extends to the graph of an everywhere defined
contraction.}\\

In the complex case, this is \cite[Theorem 1.4.1 p. 27]{adrs}. We
follow the arguments there. We consider the orthogonal projection
from $\mathcal P_2$ onto $\mathscr W_-$. Let $T$ be the densely
defined contraction with graph the relation $R$. There exist
$\mathbb H$-valued right linear functionals $c_1,\ldots,
c_\kappa$, defined on the domain of $R$, and such that
\[
Tv=\sum_{n=1}^\kappa f_nc_n(v)+w_+,
\]
where $w_+\in\mathscr W_+$ is such that $[f_n,w_+]_2=0$ for
$n=1,2,\ldots,\kappa$. Assume that $c_1$ is not bounded on its
domain, let $v_+$ be such that $c_1(v_+)=1$, and let $v_n$ be
vectors in $\mathscr V_+$ such that $c_1(v_n)=1$ and
$\lim_{n\rightarrow\infty}[v_+-v_n,v_+-v_n]_1=0$. Then $v_+$
belongs to the closure of $\ker c_1$ and so, we have that the
closure of $\ker c_1=\mathscr V_+$. Thus $\ker c_1$ contains a
strictly negative subspace of dimension $\kappa$, say $\mathscr
K_-$. For $v\in\mathscr K_-$, we have
\[
Tv=\sum_{n=2}^nf_nc_n(v).
\]
This contradicts STEP 2 and thus completes the proof of the theorem.
\end{proof}

\section{Negative squares}
\setcounter{equation}{0}
\label{sec:neqsq}
The notion of kernels with a finite number of negative squares
extend the notion of positive definite kernels. For this notion
in the quaternionic case, we send the reader to \cite[\S11]{as3}.
We recall that a $\mathbb H^{N\times N}$ Hermitian matrix $A$ has
only real (right) eigenvalues. We denote by ${\rm sq}_-(A)$ the number of
its strictly negative eigenvalues (if any).
\begin{Dn}
Let $\kappa\in\mathbb N_0$. A $\mathbb H^{N\times N}$-valued
function $K(z,\w)$ defined in a set $\Omega$ is said to have
$\kappa$ negative squares if it is Hermitian:
\begin{equation}
\label{def:herm}
K(z,w)=K(w,z)^*,\quad\forall z,w\in\Omega,
\end{equation}
and if, for every $N\in\mathbb N$ and every choice of
$z_1,\ldots, z_N\in\Omega$ and $c_1,\ldots, c_N\in\mathbb H^N$,
the $N\times N$ Hermitian matrix with $\ell, j$ entry equal to
$c_\ell^*K(z_\ell, z_j)c_j$, has at most $\kappa$ strictly
negative eigenvalues, and has exactly $\kappa$ strictly negative
eigenvalues for some choice of $N,z_1,\ldots, z_N,c_1,\ldots,
c_N$.
\end{Dn}

We will usually use the term {\sl kernel} rather than {\sl
function} to denote such $K(z,w)$'s. When $\kappa=0$, the kernel
$K(z,w)$ is positive definite. In the following theorem we recall
the following quaternionic counterparts of results well known in
the complex case. First a definition. A positive definite
$\mathbb H^{N\times N}$-valued function $Q(z,w)$ is said to be of
{\sl finite rank} if it can be factored as
\[
Q(z,w)=N(z,w)^*N(z,w),
\]
where $N(z,w)$ is $\mathbb H^{N\times M}$-valued for some
$M\in\mathbb N$. The smallest such $M$ is called the rank of $Q$.
\begin{Tm}\mbox{}\\
$(a)$ Let $K(z,w)$ be an Hermitian $\mathbb H^{N\times N}$-valued
function (see \eqref{def:herm}) for $z,w$ in some set $\Omega$.
Then, $K$ has $\kappa$ negative squares if and only if it can be
written as a difference
\[
K(z,w)=K_+(z,w)-K_-(z,w),
\]
where both $K_+$ and $K_-$ are positive definite in $\Omega$,
with moreover $K_-$ of finite rank.\\
$(b)$ There is a one-to-one correspondence between right
quaternionic reproducing kernel Pontryagin spaces of index
$\kappa$, of $\mathbb H^N$-valued functions on a set $\Omega$,
and $\mathbb H^{N\times N}$-valued functions with $\kappa$
negative squares in $\Omega$.
\end{Tm}

For a proof of these facts, see \cite[Theorems 11.5, p. 466 and
13.1, p. 472]{as3}.

\begin{Tm}
Let $K(p,q)$ be a  $\mathbb{H}^{N\times N}$-valued function with
$\kappa$ negative squares in an open nonempty subset $\Omega$ of
$\mathbb{H}$. Then there exists a unique right quaternionic
reproducing kernel Pontryagin  space $\mathscr{P}$ consisting of
$\mathbb{H}^{N}$-valued function slice hyperholomorphic in
$\Omega$ and with reproducing kernel $K(p,q)$.
\end{Tm}
\begin{proof}
The fact that there exists a unique  Pontryagin  space
$\mathscr{P}$ associated to $K$ follows as in Theorem 13.1 in
\cite{as3}. We have to show that the elements in $\mathscr{P}$
are slice hyperholomorphic. Let $\stackrel{\circ}{\mathscr{P}(K)}$ be the
linear span of the functions of the form $p\mapsto K(p,q)a$ where
$q\in \Omega$ and $a\in \mathbb{H}^N$. Since $K$ has $\kappa$
negative squares, $\stackrel{\circ}{\mathscr{P}(K)}$ has a maximal strictly
negative subspace $\mathscr{N}_{-}$ of dimension $\kappa$. By
Proposition 10.3 in \cite{as3} it is possible to write
$$
\stackrel{\circ}{\mathscr{P}(K)}=\mathscr{N}_{-} + \mathscr{N}_{-}^{[\perp]},
$$
where $\mathscr{N}_{-}^{[\perp]}$ is a quaternionic pre-Hilbert
space. The space $\mathscr{N}_{-}^{[\perp]}$ has a unique
completion, denoted by $\mathscr{N}_{+}$. Let us define
$$
\mathscr{P}:= \mathscr{N}_{+}+ \mathscr{N}_{-},
$$
with the inner product
$$
[f,f]:=[f_+,f_+]_{\mathscr{N}_{+}}+[f_-,f_-]_{\mathscr{N}_{-}}, \
\ \ \ {\rm where}\ f=f_++f_-, \  \ f_\pm \in  \mathscr{N}_{\pm}.
$$
If $f_1,\ldots, f_\kappa$ is an orthonormal basis of
$\mathscr{N}_{-}$, then
\begin{equation}\label{kerN-perp}
K(p,q)-\sum_{j=1}^\kappa f_j(p)f_j(q)^*
\end{equation}
is a reproducing kernel for $\mathscr{N}_{+}$. The functions
$f_j(p)$ are clearly slice hyperholomorphic in $p$ since they
belong to $\stackrel{\circ}{\mathscr{P}(K)}$ and so are the products
$f_j(p)f_j(q)^*$ as well as the kernel (\ref{kerN-perp}).
Therefore the elements in $\mathscr{N}_{+}$ are slice
hyperholomorphic and so are the elements in $\mathscr{P}$. This
concludes the proof.

\end{proof}

\section{Generalized Schur functions}
\setcounter{equation}{0}
\begin{Dn}
Let $J_1$ and $J_2$ be two signature matrices, respectively in
$\mathbb H^{N\times N}$ and $\mathbb H^{M\times M}$, and assume
that ${\rm sq}_- J_1={\rm sq}_-J_2$. A $\mathbb H^{N\times M}$-valued function $\Theta$, slice hyperholomorphic in a neighborhood
$\mathcal{V}$ of the origin, is called a generalized Schur function if the
kernel
\[
K_\Theta(p,q)=\sum_{\ell=0}^\infty
p^\ell(J_2-\Theta(p)J_1\Theta(q)^*)\overline{q}^\ell
\]
has a finite number, say $\kappa$, of negative squares in $\mathcal{V}$.
\end{Dn}
We will use the notation  $\mathscr S_\kappa(J_2,J_1)$ for the
class of such functions. When $N=M=1$, $\kappa=0$ and $J_1=J_2=1$,
this class was introduced in \cite{acs1}. In the statement, a
pair of operators $(C,A)$ between appropriate spaces is called
observable if
\begin{equation}
\label{obser}
\cap_{n=0}^\infty \ker CA^n=\left\{0\right\}.
\end{equation}
\begin{Tm}
Let $\Theta$ be slice hyperholomorphic in a neighborhood of the origin.
Then, it is in $\mathscr S_\kappa(J_2,J_1)$ if and only if it can
written in the form
\[
\Theta(p)=D+p C\star(I_{\mathscr P}-pA)^{-\star}B,
\]
where $\mathscr P$ is a right quaternionic Pontryagin space of
index $\kappa$, where the pair $(C,A)$ is observable, and the
operator matrix satisfies
\begin{equation}
\begin{pmatrix}A&B\\ C&D\end{pmatrix}\begin{pmatrix} I_{\mathscr P}&0\\
0&J_2\end{pmatrix}\begin{pmatrix}A&B\\ C&D\end{pmatrix}^*=
\begin{pmatrix} I_{\mathscr P}&0\\
0&J_1\end{pmatrix}.
\end{equation}
\label{tm:real}
\end{Tm}

\begin{proof} We denote by $\mathscr P(\Theta)$ the right quaternionic
reproducing kernel Pontryagin space with reproducing
$K_\Theta(p,q)$. We follow the proof of \cite[Theorem 2.2.1, p.
49]{adrs}, and we use the same densely defined linear relation as
\cite{acs1}, but this time in $({\mathscr
P}(\Theta)\oplus \mathbb H^M_{J_2})\times ({\mathscr
P}(\Theta)\oplus \mathbb H^N_{J_1})$. More precisely, $R$ is now
\[
\left\{ \begin{pmatrix} K_\Theta(p,q)\overline{q}u\\
\overline{q}v\end{pmatrix},
\begin{pmatrix}(K_\Theta(p,q)-K_\Theta(p,0))u+K_\Theta(p,0)\overline{q}v\\
(\Theta(q)^*-\Theta(0)^*)u+
\Theta(0)^*\overline{q}v\end{pmatrix}\right\}.
\]
Since ${\rm sq}_-(J_1)={\rm sq}_-(J_2)$, these Pontryagin spaces
have same negative index, and we then use Shmulyan's result to
conclude. The arguments are similar to those in \cite{acs1} and
will be omitted.
\end{proof}

We now characterize finite dimensional $\mathscr P(s)$ spaces. We
begin with a preliminary proposition.

\begin{Pn}
\label{ag}
Let
\begin{equation}
\label{relation}
\begin{pmatrix}A&B\\ C&D\end{pmatrix}
\begin{pmatrix}H&0\\ 0&J_1\end{pmatrix}
\begin{pmatrix}A&B\\ C&D\end{pmatrix}^*=
\begin{pmatrix}H&0\\ 0&J_2\end{pmatrix}
\end{equation}
and
\begin{equation}
\label{eq:s}
s(p)=D+pC\star(I-pA)^{-\star}B.
\end{equation}
Then it holds that
\[
J_2- s(p)J_1s(q)^*= C\star (I-pA)^{-\star}\star( H-p H\bar{q})\star_r (I-qA)^{-\star}\star_r C^*.
\]
\end{Pn}
\begin{proof} We rewrite the matrix identity \eqref{relation} as:
\[
\begin{split}
&J_2-DJ_1D^*=CHC^*\\
& BJ_1B^*=H-AHA^*\\
&AHC^*=-BJ_1D^*.\\
\end{split}
\]
In the sequel, for the sake of simplicity, we will write $\star$
instead of $\star_l$. Let $s(p)$ be given by \eqref{eq:s}, and
consider the function $J_2- s(p)J_1s(q)^*$ which is slice
hyperholomorphic in $p$ and $\overline{q}$ on the left and on the
right, respectively. Let us compute
\[
\begin{split}
& J_2- s(p)J_1s(q)^*=\\
& =J_2-(D+pC\star (I-pA)^{-\star}B) J_1 (D+qC\star (I-qA)^{-\star}B)^*.\\
\end{split}
\]
In order to preserve the hyperholomorphicity in $p$,
$\overline{q}$ we take, accordingly, the $\star$- product in $p$
and $\star_r$-product in $\overline{q}$ and we obtain:
\[
\begin{split}
& J_2- s(p)J_1s(q)^*\\
& =J_2-(D+pC\star (I-pA)^{-\star}B) J_1 (D^*+B^*\star_r ((I-qA)^{-\star})^*
\star_r C^* \bar{q})\\
& =J_2-DJ_1D^*- pC\star (I-pA)^{-\star} BJD^* -DJ_1B^*
\star_r((I-qA)^{-\star})^*\star_r C^*\bar{q} \\
&-pC\star (I-pA)^{-\star} BJ_1B^* \star_r((I-qA)^{-\star})^* \star_r C^*\bar{q}.\\
\end{split}
\]
Using the relations implied by (\ref{relation}) and the
identities (\ref{identities}), we obtain
\[
\begin{split}
& J_2- s(p)J_1s(q)^*\\
&=CHC^*+ pC \star (I-pA)^{-\star} AHC^* +CHA^* \star_r ((I-qA)^{-\star})^*
\star_r C^* \bar{q}\\
&-pC\star (I-pA)^{-\star} (H-AHA^*)\star_r ((I-qA)^{-\star})^* \star_r C^*\bar{q}\\
&=C\star (I-pA)^{-\star}\star\Big[ (I-pA)HC^*+ pAHC^*+
(I-pA)HA^*\star_r((I-qA)^{-\star})^* \star_r C^*\bar{q}\\
& -p(H-AHA^*)\star_r((I-qA)^{-\star})^*\star_r C^*\bar{q}\Big]\\
&=C\star (I-pA)^{-\star}\star\Big[(I-pA)H(I-qA)^*+ pAH (I-qA)^*+(I-pA)HA^*\bar{q} \\
&-p(H-AHA^*)\bar{q}\Big]\star_r ((I-qA)^{-\star})^*\star_r C^*\\
&=C\star (I-pA)^{-\star}\star\Big[ H-HA^*\bar{q}-pAH+pAHA^*\bar{q}+pAH\\
&-pAHA^*\bar{q}+HA^*\bar{q}-pAHA^*\bar{q}-pH\bar{q}+pAHA^*\bar{q}
\Big]\star_r ((I-qA)^{-\star})^*\star_r C^*\\
&=C\star (I-pA)^{-\star}\star (H-pH\bar{q})\star_r
((I-qA)^{-\star})^*\star_r C^*.
\end{split}
\]
We can also write, in an equivalent way:
\[
\begin{split}
& J_2- s(p)J_1s(q)^*\\
&=C\star (I-pA)^{-\star} H\star (1-p\bar{q})\star_r (I-qA)^{-\star}\star_r C^*\\
&=C\star (I-pA)^{-\star}\star (1-p\bar{q})\star_r H ((I-qA)^{-\star})^*\star_r C^*,\\
\end{split}
\]
or
\[
\begin{split}
& J_2- s(p)J_1s(q)^*\\
& =(C\star (I-pA)^{-\star})\star (H-pH\bar{q})\star_r (C\star
(I-qA)^{-\star})^*.
\end{split}
\]
\end{proof}

Specializing Theorem \ref{tm:real} to the finite dimensional case
we obtain:

\begin{Tm}
Let $s$ be a generalized Schur function. The associated space
right reproducing kernel Pontryagin space $\mathscr P(s)$ is
finite dimensional if and only there exists a finite dimensional
right Pontryagin space $\mathscr P$ such that:
\[
s(p)=D+pC\star (I-pA)^{-\star}B,
\]
where
\[
\begin{pmatrix}A&B\\ C&D\end{pmatrix}\,\, :\,\,
\mathscr P\oplus \mathbb H^M_{J_2}\longrightarrow \mathscr P\oplus
\mathbb H^N_{J_1}
\]
is coisometric, that is:
\begin{equation}
\begin{pmatrix}A&B\\ C&D\end{pmatrix}\begin{pmatrix} I_{\mathscr P}&0\\
0&J_2\end{pmatrix}\begin{pmatrix}A&B\\ C&D\end{pmatrix}^*=
\begin{pmatrix} I_{\mathscr P}&0\\
0&J_1\end{pmatrix}.
\end{equation}
\end{Tm}

\begin{proof} One half of the theorem follows from the preceding
proposition, while the other half is a special case of Theorem
\ref{tm:real}.

\end{proof}

%
Here we focus on the case $M=N$ and $\mathscr P(s)$ finite
dimensional.
\begin{Dn}
Let $J\in\mathbb H^{N\times N}$ be a signature function. The
$\mathbb H^{N\times N}$-valued generalized function $s$ belongs to
$s\in\mathcal U_\kappa(J)$ if the space $\mathscr P(s)$ is finite
dimensional and if ${\rm sq}_-(s)=\kappa$.
\end{Dn}

\begin{Tm}
$s\in\mathcal U_\kappa(J)$  and it is slice hyperholomorphic in a
neighborhood of the origin and only if it admits a realization
\[
s(p)=D+pC\star(I-pA)^{-\star}B
\]
where $A,B,C$ and $D$ are matrices such that
\[
\begin{pmatrix}
A&B\\C&D\end{pmatrix}\begin{pmatrix}H&0\\0&J\end{pmatrix}
\begin{pmatrix}
A&B\\C&D\end{pmatrix}^*=\begin{pmatrix}H&0\\0&J\end{pmatrix}
\]
for some Hermitian matrix $H\in\mathbb H^{N\times N}$.
\end{Tm}

\begin{proof}
First of all we observe that, if $f(p)$ is a (left) slice
hyperholomorphic function, and $C$ is a matrix, we have the
following identities which immediately follow from the definition
of (left) slice hyperholomorphic product:
\begin{equation}\label{identities}
pC\star_l f(p)= C\star_l pf(p)=C\star_l f(p)\star_l p.
\end{equation}
An analogous property holds for the right slice hyperholomorphic
product. It is also useful to recall that (compare with Section
3), if $f,g$ are left slice hyperholomorphic functions, then
$(f\star_l g)^*=g^*\star_r f^*$ and that $f\star_l C =fC$, and
analogously, if $h$ is right slice hyperholomorphic then
$C\star_r h=Ch$.
\end{proof}

For the complex-valued counterparts of the results in this
section we refer to \cite{ag,ad3}. These  last papers also
suggest factorization results, which will be considered elsewhere.

\section{Generalized Carath\'eodory functions}
\setcounter{equation}{0}
To conclude this paper we briefly study
the counterparts of the kernels \eqref{eq:kerneltheta1}.
\begin{Dn}
Let $J\in\mathbb H^{N\times N}$ be a signature matrix. A $\mathbb
C^{N\times N}$-valued function $\varphi$ slice hyperholomorphic
in a neighborhood $\mathcal{V}$ of the origin is called a generalized
Carath\'eodory function if the kernel
\[
k_\varphi(p,q)=\sum_{\ell=0}^\infty
p^\ell(\varphi(p)J+J\varphi(q)^*)\overline{q}^\ell
\]
has a finite number, say $\kappa$, of negative squares in $\mathcal{V}$.
\end{Dn}
We will use the notation  $\mathscr C_\kappa(J)$ for the class of
such functions. In the case of analytic functions, and for $N=1$
and $\kappa=0$, these functions appear in particular in the work
of Herglotz, see \cite{herglotz}, \cite{hspnw}. Still for
analytic functions, these classes were introduced and studied by
Krein and Langer, also in the operator-valued case. See
\cite{kl1}. We now give a realization theorem for such functions,
which is the counterpart in the present setting of a result of
Krein and Langer, see \cite{kl1}. As for the realization of
generalized Schur functions, we build a densely defined relation,
and apply Shmulyan's theorem (Theorem \ref{schmu} above). We
follow in the present setting the arguments in \cite[Theorem 5.2,
p. 708]{atv1}. For the notion of observability in the statement of
the theorem, see \eqref{obser}. It is equivalent to:

\begin{equation}
\label{observability}
C\star(I-pV)^{-\star}f\equiv
0\Longrightarrow f=0.
\end{equation}

\begin{Tm}
A $\mathbb C^{N\times N}$-valued function $\varphi$ slice
hyperholomorphic in a neighborhood $\mathcal V$ of the origin belongs to
$\mathscr C_\kappa(J)$ if and only if it can be written as
\begin{equation}
\label{phi}
\varphi(p)=\frac{1}{2}C\star (I_{\mathscr P}+pV)\star (I_{\mathscr P}-pV)^{-\star}C^*J
+\frac{\varphi(0)-J\varphi(0)^*J}{2},
\end{equation}
where $\mathscr P$ is a right quaternionic Pontryagin space of
index $\kappa$, $V$ is a co-isometry in $\mathscr P$, and $C$ is
a bounded operator from $\mathscr P$ into $\mathbb H^N$, and the
pair $(C,A)$ is observable.
\end{Tm}

\begin{proof}
Let $\mathscr L(\varphi)$ denote the reproducing kernel
right quaternionic Pontryagin space of functions slice
hyperholomorphic in $\mathcal V$, with reproducing kernel
$k_\varphi(p,q)$, and proceed in a number of steps. Note that in the sequel, for the sake of simplicity,
 we will write $I$ to denote the identity without specifying the space on which it is defined.\\

STEP 1: {\sl The linear relation consisting of the pairs $(F,G)\in
\mathscr L(\varphi)\times \mathscr L(\varphi)$ with
\[
F(p)=\sum_{j=1}^nk_\varphi(p,p_j)\overline{p_j}b_j,\quad{and}\quad
G(p)= \sum_{j=1}^nk_\varphi(p,p_j)b_j-k_\varphi(p,0)\left(
\sum_{\ell=1}^nb_\ell\right),
\]
where $n$ varies in $\mathbb N$, $p_1,\ldots , p_n\in \mathcal V\subset
\mathbb H$ and $b_1,\ldots, b_n\in\mathbb H^N$ is isometric, and where
by $\overline{p_j}b_j$ we mean multiplication on the right by $p_j$ on all
the components of $b_j$.}\\

We need to check that
\begin{equation}
\label{iso}
[F,F]_{\mathscr L(\varphi)}=[G,G]_{\mathscr
L(\varphi)}.
\end{equation}
We have
\[
\begin{split}
[F,F]_{\mathscr
L(\varphi)}&=[\sum_{j=1}^nk_\varphi(p,p_j)\overline{p_j}b_j,
\sum_{k=1}^nk_\varphi(p,p_k)\overline{p_k}b_k]_{\mathscr
L(\varphi)}\\
&= \sum_{j,k=1}^nb_k^*p_kk_\varphi(p_k,p_j)\overline{p_j}b_j\\
&=\sum_{\ell=1}^\infty
\sum_{j,k=1}^nb_k^*p_k^{\ell+1}(\varphi(p_k)J+J\varphi(p_j)^*)
\overline{p_j}^{\ell+1}b_j,
\end{split}
\]
while the inner product $[G,G]_{\mathscr L(\varphi)}$ is a sum of
four terms: The first is
\[
\sum_{j,k=1}^nb_k^*k_\varphi(p_k,p_j)b_j =\sum_{\ell=1}^\infty
\sum_{j,k=1}^nb_k^*p_k^{\ell}(\varphi(p_k)J+J\varphi(p_j)^*)
\overline{p_j}^{\ell}b_j.
\]
Let
\[
b=\sum_{\ell=1}^nb_\ell.
\]
The second and third terms are
\[
\begin{split}
-\left(\sum_{k=1}^nb_k^*k_\varphi(p_k,0)\right)b&=-
\sum_{k=1}^nb_k^*(\varphi(p_k)J+J\varphi(0)^*)b\\
&= - \left(\sum_{k=1}^nb_k^*\varphi(p_k)J\right)b-b^*J\varphi(0)^*b,
\end{split}
\]
and
\[
\begin{split}
-b^*\left(\sum_{k=1}^nk_\varphi(0,p_j)b_j\right)&=
\sum_{k=1}^nb^*(\varphi(0)J+J\varphi(p_j)^*)b_k\\
&=
-b^*\varphi(0)Jb- b^*\left(\sum_{j=1}^nJ\varphi(p_j)^* b_j\right),
\end{split}
\]
respectively, and the fourth term is
\[
b^*k_\varphi(0,0)b= b^*(\varphi(0)J+J\varphi(0)^*)b.
\]
Equation \eqref{iso} follows since
\[
\begin{split}
[F,F]_{\mathscr L(\varphi)}-
\sum_{j,k=1}^nb_k^*k_\varphi(p_k,p_j)b_j &=\sum_{\ell=1}^\infty
\sum_{j,k=1}^nb_k^*p_k^{\ell+1}(\varphi(p_k)J+J\varphi(p_j)^*)
\overline{p_j}^{\ell+1}b_j-\\
&\hspace{5mm}-\sum_{\ell=1}^\infty
\sum_{j,k=1}^nb_k^*p_k^{\ell}(\varphi(p_k)J+J\varphi(p_j)^*)
\overline{p_j}^{\ell}b_j\\
&=\sum_{j,k=1}^nb_k^*(\varphi(p_k)J+J\varphi(p_j)^*)b_j.
\end{split}
\]

The domain of $R$ is dense. Thus by Shmulyan's theorem (Theorem
\ref{schmu} above), $R$ is the graph of a densely defined
isometry, which extends to an isometry to all of $\mathscr
L(\varphi)$. We denote by $T$
this extension.\\

STEP 2: {\sl We compute the adjoint of the operator $T$.}\\

Let $f\in \mathscr L(\varphi)$, $h\in\mathbb H^N$ and $p\in \mathcal V$. We have:
\[
\begin{split}
h^*{p}\left((T^*f)(p)\right)&=[T^*f,k_\varphi(\cdot,
p)\overline{p}h]_{\mathscr L(\varphi)}\\&=[f\,
,\,T(k_\varphi(\cdot, p)h)]_{\mathscr
L(\varphi)}\\
&=[f\, ,\, k_\varphi(\cdot, p)h-k_\varphi(\cdot,0)h]_{\mathscr
L(\varphi)}\\
&=h^*\left(f(p)-f(0)\right),
\end{split}
\]
and hence (with $f(p)=\sum_{\ell=0}^\infty p^\ell f_\ell$)
\[
(T^*f)(p)=\begin{cases}p^{-1}(f(p)-f(0)),\quad p\not=0,\\
f_1,\quad\hspace{2.6cm} p=0.\end{cases}
\]

STEP 3: {\sl Formula \eqref{phi} is in force.}\\

We first note that $f_\ell=CR_0^\ell f$, and so
\[
f(p)=\sum_{\ell=0}^\infty p^\ell CR_0^\ell f=C\star(I-pR_0)^{-\star}f.
\]
Applying this formula to the function $C^*1=k_\varphi(\cdot, 0)$ we obtain:
\[
\varphi(p)J+J\varphi(0)^*=C\star(I-pR_0)^{-\star}C^*1\quad{\rm and}\quad
\varphi(0)J+J\varphi(0)^*=CC^*1.
\]
Multiplying the second equality by $1/2$ and making the difference with the first equality we get
\[
\varphi(p)J+\frac{1}{2}(J\varphi(0)^*-\varphi(0)J)=\frac{1}{2}C\star (I-pR_0)^{-\star} \star(I+pR_0)C^*.
\]

STEP 4: {\sl We show that conversely, every function of the form \eqref{phi} is in
$\mathscr C_\kappa(J)$.}\\

From \eqref{phi} we obtain
\begin{equation}
\label{kernel1}
\varphi(p)J+J\varphi(q)^*J= C\star(I-pV)^{-\star}\star(1-p\overline{q})\star_r((I-qV)^{-\star})^*\star_rC^*.
\end{equation}
So the reproducing kernel of $\mathscr L(\varphi)$ can be written as
\[
k_\varphi(p,q)=C\star(I-pV)^{-\star}((I-qV)^{-\star})^*\star_rC^*,
\]
since, in view of \eqref{kernel1}, the right side of the above equation satisfies
\[
k_\varphi(p,q)-pk_\varphi(p,q)\overline{q}=\varphi(p)J+J\varphi(q)^*J.
\]
In view of \eqref{observability}, $\mathscr L(\varphi)$ consists of the functions of the form
\[
f(p)=C\star(I-pV)^{-\star}\xi,\quad \xi\in\mathscr P,
\]
with the inner product
\[
[f,g]_{\mathscr L(\varphi)}=[\xi,\eta]_{\mathscr P}\quad ({\rm with}\quad g(p)=C\star(I-pV)^{-\star}\eta),
\]
and so the kernel $k_\varphi$ has exactly $\kappa$ negative squares.
\end{proof}

\begin{Cy}
When $J=I_N$ and $\kappa=0$, the function $\varphi$ has a slice
holomorphic extension to all of the unit ball of $\mathbb H$.
\end{Cy}

\begin{proof}
This follows from \eqref{phi} since $V$ is then contractive.
\end{proof}

When $J=I_N$, and in the complex variable setting generalized
Carath\'eodory functions admit another representation, namely
\begin{equation}
\label{g}
\varphi(z)=g(z)\varphi_0(z)g(1/\overline{z})^*,
\end{equation}
where $\varphi_0$ is a Carath\'eodory function (that is, the
corresponding kernel is positive definite) and $g$ is analytic
and invertible in the open unit disk. See
\cite{MR904483,MR1736921,MR1771251}. We note that in the rational
case, generalized Carath\'eodory functions are called generalized
positive functions, and play an important role in linear system
theory. We refer to \cite{a_lew_1} for a survey of the literature
and a constructive proof of the factorization \eqref{g} (in the
half-line case) in the scalar rational case.

\end{document}